\documentclass[11pt, leqno]{amsart}

\usepackage[bookmarks=true,colorlinks=true,pdfstartview=FitV,linkcolor=blue, citecolor=red,urlcolor=green]{hyperref}  

\usepackage{graphicx}
\usepackage{amssymb,amsmath,mathabx}
\usepackage{epstopdf}
\usepackage{url}
\usepackage{bm}

\usepackage{lpic}

\usepackage{subfigure}

\usepackage{syntonly}

\usepackage[usenames,dvipsnames,svgnames,table]{xcolor}

\newcommand{\deltxt}[1]{}

\newtheorem{theorem}{Theorem}[section]
\newtheorem{proposition}[theorem]{Proposition}
\newtheorem{lemma}[theorem]{Lemma}
\newtheorem{corollary}[theorem]{Corollary}

\numberwithin{equation}{section}

\theoremstyle{remark}
\newtheorem{rem}[theorem]{Remark}
\newtheorem{definition}[theorem]{Definition}
\newtheorem{hyp}[theorem]{Hypothesis}

\newcommand\bR{{\mathbb{R}}}
\newcommand\bC{{\mathbb C}}
\newcommand\bZ{{\mathbb Z}}

\newcommand\bN{{\mathbb N}}

\newcommand\bH{{\mathbb{H}}}

\newcommand\dev{{\bf dev}}
\newcommand\SI{{\mathbb{S}}}

\newcommand\inte{{\rm int}}
\newcommand\Bd{{\rm bd}}
\newcommand\clo{{\rm Cl}}
\newcommand\bdd{{\mathbf{d}}}

\newcommand\ra{\rightarrow}

\newcommand\che{\check}
\newcommand\emp{\emptyset}
\newcommand\eps{\epsilon}

\newcommand\Aff{{\mathbf{Aff}}}

\newcommand\Idd{{\rm I}}

\newcommand\SL{{\mathsf{SL}}}
\newcommand\PSL{{\mathsf{PSL}}}

\newcommand\PGL{{\mathsf{PGL}}}

\newcommand\GL{{\mathsf{GL}}}

\newcommand\Mod{{\mathrm{Mod}}}

\newcommand{\llrrparen}[1]{
  \left(\mkern-4mu\left(#1\right)\mkern-4mu\right)}

\newcommand{\bemph}[1]{{\upshape#1}} 
\newcommand{\ep}[1]{\bemph{(}#1\,\bemph{)}} 

\begin{document}

\title{Convex and concave decompositions of affine $3$-manifolds}


\author{Suhyoung Choi} 
\address{ Department of Mathematics \\ KAIST \\
	Yuseong-gu 291 Daehak-ro \\
	Daejeon 34141, South Korea \\
}
\email{schoi@math.kaist.ac.kr}

\thanks{This work was supported by the National Research Foundation
	of Korea (NRF) grant funded by the Korea government (MEST) (No.2013R1A1A2056698)}
















\date{\today}

 
\begin{abstract} 
A (flat) affine $3$-manifold is a $3$-manifold with an atlas of charts to an affine space 
$\bR^3$ with transition maps in the affine transformation group $\Aff(\bR^3)$.
We will show that a connected closed affine $3$-manifold is either an affine Hopf $3$-manifold or
decomposes canonically to concave affine submanifolds with incompressible boundary,
toral $\pi$-submanifolds and $2$-convex affine manifolds, each of which
is an irreducible $3$-manifold.  It follows that 
if there is no toral $\pi$-submanifold, then $M$ is prime. 
Finally, we prove
that if a closed affine manifold is covered by a domain in $\bR^{3}$, 
then $M$ is irreducible or is an affine Hopf manifold. 

%

\end{abstract}



%


\maketitle

\section{Introduction}

\subsection{Introduction and history} 
An \hyperlink{term-afm}{affine manifold} is a manifold with an atlas of charts to $\bR^{n}$, $n \geq 2$, where the transition maps are in the affine group. 
Euclidean manifolds are examples. A \hyperlink{term-hopf}{Hopf manifold} that is the quotient of $\bR^{n} -\{O\}$ by a linear contraction group, i.e., 
a group of linear transformation generated by an element with eigenvalues of norm greater than $1$
is an example. A half-Hopf manifold is the quotient of $U-\{O\}$ by a linear contraction group  for a closed upper half-space $U$ of $\bR^{n}$. 
(See Proposition \ref{prop:Hopf}.)  



For the currently most extensive set of examples of affine manifolds, see the paper by Sullivan and Thurston \cite{SuTh}. 
We still have not obtained essentially different examples to theirs to this date. (See also Carri\`ere \cite{carr}, 
Smillie \cite{Smillie} and Benoist \cite{BenTor} and \cite{BenNil}.) 
A connected compact affine $3$-manifold is {\em radiant} if the holonomy group fixes a unique point.
(See Section \ref{sub:ram} and  Barbot \cite{Barbot}, 
 Fried, Goldman, and Hirsch \cite{FGH}.) 
Such a manifold has a complete flow called a {\em radiant flow}. 
A generalized affine suspension is a radiant affine manifold admitting a total cross-section. (See Proposition \ref{prop:tcs}.)
A radiant affine $n$-manifold can be constructed easily from a real projective $(n-1)$-manifold using generalized affine suspension. 
(See Section 2.2 of \cite{Barbot} or Chapter 3 of \cite{rdsv}.)

A $3$-manifold $M$ is {\em prime} if $M$ is a connected sum of two manifolds $M_{1}$ and $M_{2}$, then
$M_{1}$ or $M_{2}$ is homeomorphic to a $3$-sphere. 
The subject of this paper is the following: 
The question of Goldman in Problem 6 in the Open problems section of \cite{GTP} is 
whether closed affine $3$-manifolds are prime. 
We showed that $2$-convex affine $3$-manifolds are irreducible in \cite{uaf}. 
Our Theorem \ref{thm:main3} shows that closed affine manifolds may 
be obtainable by gluing \hyperlink{term-tsp}{toral $\pi$-submanifolds} 
which are solid tori or solid Klein bottles with special geometric properties 
to irreducible $3$-manifolds.
This construction may result in reducible $3$-manifolds as we can see from Gordon \cite{Go}.
Hence, the nonexistence of  solid tori or solid Klein bottles with special geometric properties in 
a closed affine $3$-manifold $M$ would show that $M$ is prime. 
(See Corollary \ref{cor:main4}.) 
We question whether toral $\pi$-submanifolds can occur at all.
We also answer the question when $M$ is covered by a domain in an affine space 
by Corollary \ref{cor:main5}. 


For the related real projective structures on closed $3$-manifolds, Cooper and Goldman \cite{CG} showed that 
 a connected sum $\bR P^3 \hash \bR P^3$ admits no real projective structure. 
 For these topics, a good reference is given by Goldman \cite{Gnote}, 
 originally given as lecture notes in the 1980s.

\subsection{Main results} 
We give some definitions which we will give more precisely later. 
A \hyperlink{term-rpm}{real projective structure} on a manifold $M$ is 
a maximal atlas of charts to $\bR P^{n}$ with transition maps in the projective group $\PGL(n+1, \bR)$. 
$M$ is called a {\em real projective manifold}. 

We use the double-covering map $\SI^{n} \ra \bR P^n$, and hence $\SI^{n}$ has a real projective structure. 
The group of projective automorphism of $\bR P^{n}$ is $\PGL(n+1, \bR)$ and that of $\SI^{n}$ is $\SL_{\pm}(n+1, \bR)$. 

We recall the main results of \cite{psconv} which we will state in Section \ref{subsec:convconc} in a more detailed way. 
Let $M$ be a closed real projective manifold.
Let $\tilde M$ be the universal cover and $\pi_{1}(M)$ the deck transformation group. 
A real projective structure on $M$ gives us an immersion $\dev:\tilde M \ra \SI^{n}$ 
equivariant with respect to 
a homomorphism $h: \pi_{1}(M) \ra \SL_{\pm}(n+1, \bR)$.  The real projective structure gives these data.

Recall a group $\Aff(\bR^n)$ of affine transformations of form $x \mapsto Mx + b$ for 
$M \in \GL(n, \bR)$ and $b \in \bR^n$. 
The real projective space $\bR P^{n}$ contains the affine space $\bR^{n}$ as  a complement of a hyperspace, 
and affine transformation groups  naturally extend to projective automorphisms.  
Affine geodesics also extend to projective geodesics.

We will look an affine manifold as a \hyperlink{term-rpm}{real projective manifold}, 
i.e., a manifold with an atlas of charts to $\bR^{n}\subset \SI^{n}$ with transition maps 
in the affine group $\Aff(\bR^{n}) \subset \SL_{\pm}(n+1, \bR)$.  An affine manifold has a canonical real projective structure
since the charts and the transition maps are projective also. (The converse is not true.)

Let $K_{h}$ be the kernel of $h$, normal in $\pi_{1}(M)$. 
We cover $M$ by the holonomy cover $M_{h}= \tilde M/K_{h}$ corresponding to $K_{h}$ with
\begin{itemize}
\item  an induced and lifted immersion $\dev_{h}:M_{h} \ra \SI^{n}$ and 
\item an induced holonomy homomorphism $h_{h}: \pi_{1}(M) /K_{h}\ra \SL_{\pm}(n+1, \bR)$ satisfying 
\[ \dev_{h} \circ g = h_{h}(g) \circ \dev_{h} \hbox{ for } g \in \pi_{1}(M)/K_{h}.\]
\end{itemize} 
Let $M_{h}$ have the path metric of the Riemannian metric pulled  back from 
the Fubini-Study Riemannian metric of $\SI^{3}$.
The Cauchy completion $\che M_{h}$ of $M_{h}$ is called 
a \hyperlink{term-kcp}{{\em Kuiper completion}}.  
The \hyperlink{term-kcp}{ideal set} is $M_{h, \infty} := \che M_h - M_h$.

A \hypertarget{term-shem}{$3$-hemisphere} is a closed $3$-hemisphere in $\SI^3$, and 
a \hypertarget{term-sbih}{$3$-bihedron} is the closure of a component $H-\SI^2$ for a $3$-hemisphere $H$ 
with a great $2$-sphere $\SI^2$ passing $H^o$. These have real projective structures induced from the double-covering map $\SI^3 \ra \bR P^3$.

If the universal cover $\tilde M$ is projectively diffeomorphic to an open hemisphere, i.e., $\bR^{n}$,  
then $M$ is called a {\em complete affine manifold}.  
If the universal cover $\tilde M$ is  projectively diffeomorphic to an open $3$-bihedron, 
we call $M$ a {\em bihedral real projective manifold}.  





A \hyperlink{term-hcr}{hemispherical $3$-crescent} is a $3$-hemisphere in $\che M_{h}$ with boundary $2$-hemisphere in the ideal set. 
A \hyperlink{term-bcr}{ bihedral $3$-crescent} is a $3$-bihedron $B$ in $\che M_h$ so that 
a boundary $2$-hemisphere is the ideal set where we assume that $\che M_{h}$ has no hemispherical $3$-crescent. 
(See Section \ref{subsub:Kuiper} for definitions and Hypothesis \ref{hyp:nohem}.)
A \hyperlink{term-camI}{concave affine $3$-manifold} is a codimension-zero connected compact submanifold of $M$ defined in \cite{psconv}.
We cover these in Section \ref{subsec:crcaf}. 
The interior of a concave affine $3$-manifold has a canonical affine structure inducing its real projective structure. 
The \hyperlink{term-tfsm}{{\em two-faced  submanifold of type I}} of a real projective $3$-manifold $M$ 
is roughly given as the totally geodesic submanifold arising 
from the intersection in $M_{h}$ of two hemispherical $3$-crescents meeting only in the boundary. 
The \hyperlink{term-tfsm}{{\em two-faced  submanifold of type II}} of a real projective $3$-manifold $M$ 
is roughly defined as the totally geodesic submanifold arising 
from the intersection in $M_{h}$ of two bihedral $3$-crescents meeting only in the boundary. 
For the precise definitions, see Section \ref{subsec:crcaf}. 

Let $T$ be a convex simplex in an affine space  $\bR^3$ with faces $F_0, F_1, F_2, $ and $F_3$. 
A {\rm real projective or affine $3$-manifold} is {\em $2$-convex} if every projective map $f: T^o \cup F_1 \cup F_2 \cup F_3 \ra M$ 
extends to $f:T \ra M$. (Y. Carri\`ere \cite{carr} first defined this concept.)

\begin{theorem}[\cite{psconv}]\label{thm:psconv} 
Suppose that $M$ is a compact real projective $3$-manifold with empty or convex boundary
that is neither complete affine nor bihedral. 
Suppose that $M$ is not $2$-convex. Then $\che M_{h}$ contains a hemispherical or bihedral $3$-crescent.
\end{theorem}

Now, we sketch the process of {\em convex-concave decomposition} in \cite{psconv} which 
we recall in Section \ref{subsec:convconc} in more details: 
\begin{itemize} 
\item Suppose that a hemispherical $3$-crescent $R \subset \che M_{h}$ exists. 
\begin{itemize} 
\item  If there is the two-faced submanifold of type I, then we can \hyperlink{term-spl}{split} $M$ along this submanifold to obtain 
$M^{s}$. If not, we let $M^{s} = M$.  Let $M^{s}_{h}$ denote the corresponding cover of $M^{s}$ obtained by splitting $M_{h}$ and taking 
a union of components,
and let $\che M^{s}_{h}$ be its Kuiper completion. 
\item  Then hemispherical $3$-crescents in $\che M^{s}_{h}$ are mutually disjoint and 
their intersection with $M^{s}_{h}$ cover 
 compact submanifolds, called {\em concave affine manifolds of type I}. 
 \item  We remove all these from $M^{s}$. Then we let the resulting compact manifold be called $M^{(1)}$. The boundary is still convex.
 \end{itemize} 
\item  Let $M^{(1)}_{h}$ denote the cover of $M^{(1)}$ obtained by removing corresponding submanifolds from $M^{s}_{h}$, and let 
 $\che M^{(1)}_{h}$ be the Kuiper completion of $M^{(1)}_{h}$. 
Suppose that there is a bihedral $3$-crescent $R\subset \che M^{(1)}_{h}$. 
\begin{itemize} 
\item If there is the two-faced submanifold of type II, then we can \hyperlink{term-spl}{split} $M^{(1)}$ along this submanifold to obtain 
$M^{(1)s}$. If not, we let $M^{(1)s} = M^{(1)}$. 
\item Let $M^{(1)s}_{h}$ denote the cover of $M^{(1)s}$ obtained from $M^{(1)}_{h}$ by splitting and taking a union of components,
and let $\che M^{(1)s}_{h}$ be the Kuiper completion. 
Then the intersection of $M^{(1)s}_{h}$ with 
the union of bihedral $3$-crescents in $\che M^{(1)s}_{h}$ covers the union of a mutually disjoint collection of 
 compact submanifolds, called {\em concave affine manifolds of type II}. 
\item We remove all these from $M^{(1)s}$. Then the resulting compact real projective manifold $M^{(2)}$ with convex boundary is $2$-convex. 
\end{itemize} 
\end{itemize} 

We will further sharpen the result in this paper. 
 A \hyperlink{term-tsp}{toral $\pi$-submanifold} is a compact radiant concave affine $3$-manifold with the virtually infinite-cyclic fundamental group covered by  a special domain in a hemisphere. 
We will later show that a toral $\pi$-submanifold is homeomorphic to a solid torus or 
a solid Klein bottle. (See Definition \ref{defn:toralpi} and Lemmas \ref{lem:toralpiI} and \ref {lem:toralpi}.) 


\begin{theorem} \label{thm:main1}
Let $M$ be a connected compact real projective $3$-manifold with empty or convex boundary
that is neither complete affine nor bihedral. 
\begin{itemize}
\item Let $M^s$ be the resulting real projective $3$-manifold after \hyperlink{term-spl}{splitting} along the
two-faced totally geodesic submanifold of type I\, {\rm (}resp. of type II\,{\rm ).}
\item Let $N$ be a compact concave affine $3$-manifold in $M^s$ with compressible boundary
of type I\, {\rm (}resp. of type II\,{\rm ).}
\end{itemize}
Then $N$ is a toral $\pi$-submanifold of type I\, {\rm (}resp. contains a unique maximal toral $\pi$-submanifold of type II\,{\rm )} or $M$ is 
an \hyperlink{term-hopf}{affine Hopf $3$-manifold}.
\end{theorem}

So far, our results are on real projective $3$-manifolds. 
Now we go over to the result specific  to affine $3$-manifolds. 

\begin{theorem}\label{thm:main3} 
Let $M$ be a connected compact affine $3$-manifold with empty or convex boundary. 
Suppose that $M$ is neither complete affine nor bihedral and is not affine Hopf $3$-manifold. 
\begin{itemize}
\item Let $M^{s}$ be the resulting real projective $3$-manifold after \hyperlink{term-spl}{splitting} along 
the two-faced totally geodesic submanifold of type I. 
\item Let $M^{(1)}$ be obtained by removing all concave affine manifolds of $M^{s}$. 
$M^{s}$ decomposes into concave affine manifolds of type I with boundary  incompressible in $M^{s}$
and toral $\pi$-submanifolds of type I. 
\item Let $M^{(1)s}$ denote the $M^{(1)}$ split along the two-faced submanifold of type II. 
\end{itemize} 
Then 
$M^{(1)s}$ decomposes into compact submanifolds as follows\,{\rm :} 
\begin{itemize}
\item a $2$-convex affine $3$-manifold with convex boundary, 
\item toral $\pi$-submanifolds of type II 
in concave affine $3$-manifolds with compressible boundary with the virtually cyclic holonomy group, or 
\item concave affine $3$-manifolds of type II with 
boundary incompressible in $M^{(1)}$.
\end{itemize}
Then the finally decomposed submanifolds from the decomposition are prime $3$-manifolds.
\end{theorem}
However, the above decomposition is not necessarily a prime decomposition. 
The following gives us a criterion for the primeness. 

\begin{corollary} \label{cor:main4} 
Let $M$ be a connected compact affine $3$-manifold with empty or convex boundary. 
Suppose that there is no projectively embedded toral $\pi$-submanifold of type I or II in $M$. 
Then $M$ is irreducible or is an affine Hopf $3$-manifold and hence is prime. 
\end{corollary} 
We question whether that the above concave affine manifolds are 
maximal in the sense of \cite{psconv}. 

The following answers Goldman's question partially. 
\begin{corollary}[Choi-Wu] \label{cor:main5}
Suppose that $M$ is a connected closed affine manifold covered by a domain $\Omega$ in $\bR^{3}$.
Then $M$ is either irreducible or is an affine Hopf $3$-manifold. 
\end{corollary} 

One significance of this paper is to see how far the techniques of \cite{uaf} and \cite{psconv} can 
be applied to solve this problem. We isolated some objects here.

\subsection{Outline} 
The main tools of this paper are from three long papers \cite{psconv}, \cite{rdsv}, and \cite{uaf}. 
We summarize the results of  \cite{rdsv} and \cite{psconv} in Section \ref{sec:prel}. 
In Section \ref{sub:ram}, we recall radiant affine $n$-manifolds and recall some results of \cite{rdsv}.  
In Section \ref{sec:Hopf}, we prove various facts about affine Hopf manifolds and half-Hopf manifolds. 
In Section \ref{subsec:convconc}, we recall the convex and concave decomposition of real projective structures.
We recall $3$-crescents and two-faced submanifolds and the decomposition theory in \cite{psconv}. 

In Section \ref{sec:concave}, Theorem \ref{thm:tfaced} claims that if the two-faced submanifold
is non-$\pi_{1}$-injective, then the manifold is finitely covered by an affine Hopf $3$-manifold.
The idea for the proof is by a so-called disk-fixed-point argument, Proposition \ref{prop:diskf}; 
that is, we can find an attracting fixed point of a deck transformation $g$
using a simple closed curve $c$ bounding a disk $D$ with $g(c) \subset D^{o}$. 
We prove Theorem \ref{thm:tfaced} in Section \ref{sub:tfaced}. 

The main technical core results are 
Theorems \ref{thm:caffI} and \ref{thm:caffII} in Section \ref{subsec:caff}. 
We show that a cover of the concave affine $3$-manifold being a union of mutually intersecting 
$3$-crescents must be mapped to a domain in a hemisphere by $\dev_{h}$, and the boundary has a unique annulus component. 
Since the fundamental group of $N$ acts on an annulus
covering its boundary properly and freely, the fundamental group is virtually infinite-cyclic by Lemma \ref{lem:actann}. 
We complete the final part of the proof in Section \ref{subsec:toralpi} where we show that
these concave affine $3$-manifolds contain toral $\pi$-submanifolds. 
We also show that a toral $\pi$-submanifold is homeomorphic to a solid torus or a solid Klein bottle. 
We prove Theorem \ref{thm:main1} at the end. 

In Section \ref{sec:toraldec}, we discuss the decomposition of $M$ into $2$-convex real projective 
$3$-manifolds with convex boundary and toral $\pi$-submanifolds. 
i.e., Theorem \ref{thm:decompose}.
We use the convex and concave decomposition theorem of \cite{psconv} 
and Theorems \ref{thm:caffI} and \ref{thm:caffII} and replacing the compact concave affine $3$-manifolds with compressible boundary with 
toral $\pi$-submanifolds.  
We prove Theorem \ref{thm:main3} and Corollaries \ref{cor:main4} and \ref{cor:main5}
lastly here.  




We thank Yves Carri\`ere, David Fried, Bill Goldman and Weiqiang Wu for fruitful discussions.

\section{Preliminary} \label{sec:prel}

\subsection{ Some $3$-manifold topology } 
Let $K$ be a manifold. 
Let $\mathrm{Diff}(K)$ be the group of diffeomorphisms of $K$ with the usual $C^r$-topology, 
$r \geq 0$, and $\mathrm{Diff}_0(K)$ the identity component of this group.
We define the mapping class group $\Mod(K)$ of a manifold $K$ to be the group 
$\mathrm{Diff}(K)/\mathrm{Diff}_0(K)$. 


Since $\Mod(\SI^2) = \bZ/2\bZ$ is a classical work of Smale \cite{Sm}, 
there exist only two homeomorphism types of $\SI^{2}$-bundle over $\SI^{1}$. 
If $M'$ is orientable, then $M'$ is homeomorphic to $\SI^{2}\times \SI^{1}$. 
If not, $M'$ is a non-orientable $\SI^{2}$-bundle over $\SI^{1}$. 
The following is well known. 
\begin{lemma} \label{lem:3mfld} 
	Suppose that $\tilde N = K \times \bR$ for a compact manifold $K$   
	covers a compact manifold $N$ as a regular cover. 
	Suppose that $\Mod(K)$ is finite. 
	Then $N $ is finitely covered by $K \times \SI^1$. 
\end{lemma} 


A $n$-manifold is {\em irreducible} if every embedded two-sided 
$n-1$-sphere bounds a $3$-ball. Also, prime $3$-manifolds are either irreducible or 
are homeomorphic to a $\SI^{2}$-bundle over $\SI^{1}$.
(See Lemma 3.13 of \cite{Hempel}.)

\begin{lemma} \label{lem:vic} 
Let $L$ be a compact $3$-manifold with the universal cover whose interior is an open cell, 
and $\pi_{1}(L)$ is virtually infinite-cyclic. 
Then $L$ is homeomorphic to a solid torus or a solid Klein-bottle. 
\end{lemma} 
\begin{proof} 
Since the interior of the universal cover of $L$ is a cell, $L$ is irreducible. 
By Theorem 5.2 of \cite{Hempel}, $L$ is  finitely covered by a solid torus.
Hence, $\partial L$ is homeomorphic to a torus or a Klein bottle. 
Since $\partial L$ is not $\pi$-injective, $\partial L$ is compressible by Dehn's lemma. 
Since $\partial L$ is compressible, we can find a disk $D$ with $\partial D \subset \partial L$. 
Since $L$ is irreducible, $L-D$ is a cell. Therefore, the conclusion follows. 
\end{proof}

\begin{lemma} \label{lem:actann} 
Suppose that a group $G$ acts on an annulus $A$ faithfully, freely, and properly discontinuously. 
Then $G$ is virtually infinite-cyclic.
\end{lemma} 
\begin{proof} 
$A/G$ is a closed surface. Let $c$ be an essential simple closed curve in $A$. 
Then there is a finite index subgroup $G'$ preserving the ends of $A$ of $G$ so that $c$ is embedded to 
a simple closed curve in $A/G'$. Then $G'$ is infinite-cyclic. 
\end{proof}

\subsection{The projective geometry of the sphere}

Let $V$ be a vector space. 
Define $P(V)$ as $V -\{0\}/\sim$ where $x \sim y$ if and only if $x = s y$ for $s \in \bR -\{0\}$. 
$\PGL(V)$ acts on this space where $\PGL(\bR^n) = \PGL(n, \bR)$. 

Recall that $\bR P^{n}  = P(\bR^{n+1})$. 
A subspace of $\bR P^n$ is the image $V-\{O\}$ of a subspace $V$ of $\bR^{n+1}$ under the projection. 
The group of projective automorphisms is $\PGL(n+1, \bR)$ acting on $\bR P^n$ in the standard manner. 
A \hypertarget{term-rpm}{{\em real projective $n$-manifold with empty or convex boundary}} is a manifold with empty or nonempty boundary 
with an atlas of charts to $\bR P^n$ with transition maps in $\PGL(n+1, \bR)$ so that 
each point of the boundary has a chart to a convex domain with boundary in 
$\bR P^n$. 
A maximal atlas is called a {\em real projective structure}. 
The boundary is {\em totally geodesic} if each boundary point has a neighborhood projectively diffeomorphic 
to an open set in a half-space of an affine space meeting the boundary. 


An \hypertarget{term-afm}{{\em affine $n$-manifold with empty or convex boundary}} 
is an $n$-manifold with smooth boundary and an atlas of charts to open subsets or 
convex domains in $\bR^n$ and the transition maps in $\Aff(\bR^{n})$. 
Since the affine transformations are projective, an affine $n$-manifold has a canonical real projective structure. 
We consider such $n$-manifolds as real projective $n$-manifolds with special structures in this paper. 
A real projective manifold projectively homeomorphic to an affine manifold is called an {\em affine manifold} in this paper.

\begin{definition}\label{defn:Hopf} 
An elementary example is an \hypertarget{term-hopf}{{\em affine Hopf $n$-manifold}} that is the quotient of $\bR^n - \{O\}$ 
by an infinite-cyclic group generated by a linear map $g$ all of whose eigenvalues have norm $> 1$
or by $\langle g, -\Idd \rangle$ for $g$ as above. The quotient is a manifold by Proposition \ref{prop:hopf}. 

If $g$ acts on an $(n-1)$-plane passing $O$, and the half-space $H$ in $\bR^{n}$ bounded by it, 
then $(H -\{O\})/\langle g \rangle$  is called a {\em half-Hopf $n$-manifold}. 
A real projective manifold projectively homeomorphic to an affine Hopf $n$-manifold or 
a half-Hopf $n$-manifold is called by the same name in this paper. 
(See \cite{Hopf} for a conformally flat version and \cite{SuTh}.) 
\end{definition} 



Let $\bR_+ := \{t| t \in \bR, t> 0\}$. 
Define $\SI(V)$ as $V -\{0\}/\sim$ where $x \sim y$ if and only if $x = sy$ for $s \in \bR_+$. 
$\SL_\pm(V)$ acts on $\SI(V)$ transitively and faithfully. 
There is a double cover $\SI(V) \ra P(V)$ with the deck transformation group 
generated by the antipodal map $\mathcal{A}:\SI(V) \ra \SI(V)$ induced from 
the linear map $V \ra V$ given by $v \ra -v $. 
We denote by $\llrrparen{v}$ the equivalence class of $v$ in $\SI(V)$. 
The {\em homogeneous coordinate system} of $\SI(\bR^n)$ is given by denoting each point by
$\llrrparen{x_1, \dots, x_n}$ for the vector $(x_1, \dots, x_n) \ne 0$. 

We denote by $\SI^n$ the space $\SI(\bR^{n+1})$. 
The real projective sphere 
$\SI^n$ has a real projective structure given by the double covering map 
to $\bR P^n$. 
The group of projective automorphisms of $\SI^n$ form $\SL_\pm(n+1, \bR)$ as obtained 
by the standard action of $\GL(n+1, \bR)$ on $\bR^{n+1}$. 

We embed $\bR^{n}$ as an open $n$-hemisphere $\bH^o$ in $\SI^n$ for a closed $n$-hemisphere $\bH$
by sending $(x_1, x_2, \dots, x)$ to $\llrrparen{1, x_1, x_2, \dots, x_n}$. We identify $\bR^n$ with $\bH^o$.
The boundary of $\bR^{n}$ is a great sphere $\SI^{n-1}_{\infty}$ given by $x_0 =0$. 
The group of projective automorphisms acting on $\bH$ equals the group $\Aff(\bR^{n})$ of affine transformations 
of $\bH^o=\bR^{n}$.  
(A good reference for all these geometric topics is the book by Berger \cite{Berger}.)

We take the universal cover $\tilde M$ of a $n$-manifold $M$. 
The existence of a real projective structure on $M$ gives us 
\begin{itemize}
\item an immersion $\dev:\tilde M \ra \bR P^n$, called a {\em developing map} and 
\item a homomorphism  $h: \pi_1(M) \ra \PGL(n+1, \bR)$, called a {\em holonomy homomorphism} 
\end{itemize}
satisfying $\dev \circ \gamma = h(\gamma) \circ \dev$ for each $\gamma \in \pi_{1}(M)$. 

By lifting $\dev$, we obtain
\begin{itemize} 
\item a well-defined immersion $\dev':\tilde M \ra \SI^n$
and 
\item a homomorphism $h': \pi_1(M) \ra \SL_\pm(n+1, \bR)$
\end{itemize} 
so that $\dev' \circ g = h'(g) \circ \dev'$ for each deck transformation $g$ of $\tilde M$. 

Let $K_h$ be the kernel of $h': \pi_1(M) \ra \SL_\pm(n+1, \bR)$.
Let $M_h := \tilde M/K_h$ be a so-called holonomy cover. Then $\dev'$ induces an immersion 
$\dev_h: M_h \ra \SI^3$. The deck transformation group $\Gamma_h$ of $p_{h}:M_h \ra M$ is
isomorphic to $\pi_1(M)/K_h$. 
We obtain 
\begin{itemize} 
\item an immersion $\dev_h: M_h \ra \SI^n$, also called a {\em developing map} and
\item  a homomorphism $h_h: \pi_1(M)/K_{h} \ra  \SL_\pm(n+1, \bR)$, also called a {\em holonomy homomorphism} 
\end{itemize} 
satisfying \[\dev_h \circ \gamma = h_h(\gamma) \circ \dev_h \hbox{ for } \gamma \in \Gamma_h.\] 

\begin{lemma} \label{lem:holc} 
Consider a cover $M'$ of $M$ with a covering map $p_{M}:M' \ra M$
with a deck transformation group $\Gamma'$. Let $p_{M'}:\tilde M \ra M'$ denote the covering map
induced by the universal covering map $\tilde M \ra M$. 
Then 
\begin{itemize}
\item given a projective immersion $\dev':M' \ra \SI^{n}$ satisfying 
$\dev'= \dev \circ p_{M'}$, there is a homomorphism $h':\Gamma' \ra \SL_\pm(n+1, \bR)$ satisfying 
$\dev' \circ \gamma = h'(\gamma) \circ \dev'$ for every $\gamma \in \Gamma'$.
\item $\dev'$ is a holonomy cover if and only if 
$p_{M}:M' \ra M$ is a regular cover and $h'|\Gamma'$  is injective. 
\end{itemize} 
\end{lemma} 
\begin{proof} Straightforward. 
\end{proof} 


\begin{lemma}\label{lem:regular} 
For any connected submanifold $N$ of $M$, let $N_h$ denote a component of its inverse image in $M_h$. 
Then $p_h|N_h: N_h \ra N$ is a holonomy covering map also and the deck transformation group equals 
the subgroup $\Gamma_{h, N_h}$  of $\Gamma_h$ acting on $N_h$. 
For the developing map, $\dev_{h, N_h} = \dev_h|N_h$ holds and for the corresponding holonomy 
homomorphism, $h_{N_h}= h_h|\Gamma_{h, N_h}$ holds. 
\end{lemma} 
\begin{proof} 
First, $\Gamma_{h, N_{h}} \ra \Gamma_{h}$ is injective. Since $h_{h}|\Gamma_{h}$ is injective, 
$h_{h}|\Gamma_{h, N_{h}}$ is injective. 
Since $\Gamma_{h, N_{h}}$ is the regular deck transformation group of $p_{h}| N_{h}$, 
we are done by Lemma \ref{lem:holc}. 
\end{proof} 


\begin{lemma} \label{lem:residual} 
The deck transformation group $\Gamma_{h}$ of $M_{h}$ is residually finite. 
So, is $\Gamma_{h, N_{h}}$ for any connected submanifold $N$ of $M$. 
\end{lemma}
\begin{proof} 
Under $h_{h}$, $\Gamma_{h}$ is mapped injectively into a linear group $\SL_{\pm}(n+1, \bR)$. 
Selberg-Malcev lemma implies the conclusion.
\end{proof}

\subsection{Radiant affine $n$-manifolds} \label{sub:ram} 

Given any affine coordinates $x_{i}, i=1,\dots, n$, of $\bR^{n}$, 
a vector field $\sum_{i=1}^{n} x_{i}\frac{\partial}{\partial x_{i}}$ is called a {\em radiant} vector field.
$O$ of the coordinate system is called the {\em origin} of the radiant vector field.
Suppose that the holonomy group of an affine $n$-manifold $M$ fixes $O$. 
Then $\dev_{h}: M_{h} \ra \bR^{n}$ is an immersion and 
the radiant vector field lifts to a vector field in $M_{h}$. 
The vector field is invariant under the deck transformations of $M_{h}$
and hence induces a vector field on $M$. 
The vector field on $M$ is also called a {\em radiant vector field}. 
(see Barbot \cite{Barbot} and Chapter 3 of \cite{rdsv}.) 
This gives us a {\em radiant} flow:  
\[ \bR\times M \ra M.\]

 
 Let $M$ be a radiant affine manifold with the holonomy group fixing a point $O$. 
 A {\em radial line} in $M_{h}$ is an arc $\alpha$ in $M_{h}$ so that 
 $\dev|\alpha$ is an embedding to a component of a complete real line $l$ with $O$ removed.  
 \begin{proposition}\label{prop:rad} 
 Let $M$ be a compact affine $n$-manifold with empty or nonempty boundary.
 Suppose that the holonomy group fixes the origin of a radiant vector field, 
 and the boundary is tangent to the radiant vector field. 
 Then $\dev_{h}(M_{h})$ misses the origin of a vector field and $M_{h}$ is foliated by radial lines. 
 \end{proposition} 
\begin{proof} 
See the proof of Proposition 2.4 of Barbot \cite{Barbot}. 
\end{proof}


Let $||\cdot||$ denote the Euclidean metric of $\bR^{n}$. 
Given a real projective $(n-1)$-manifold $\Sigma$ and a projective automorphism $\phi:\Sigma \ra \Sigma$, 
we can obtain a radiant affine $n$-manifold homeomorphic to the mapping torus 
\[\Sigma\times I/\sim \hbox{ where for every } x \in \Sigma, (x, 1) \sim (\phi(x), 0).\] 
Let $\dev:\tilde \Sigma \ra \SI^{n-1} \subset \bR^{n}$ 
be a developing map with holonomy homomorphism $h: \pi_{1}(\Sigma) \ra \SL_{\pm}(n, \bR)$. 
Then we extend $\dev$ to 
\[\dev': \tilde \Sigma \times \bR \ra \bR^{n} \hbox{ by } (x, t) \mapsto \exp(t) \dev(x).\] 
For each element $\gamma$ of $\pi_{1}(\Sigma)$, we define the action of $\pi_{1}(\Sigma)$ on $\Sigma \times \bR$ by 
\[\gamma(x, t) = (\gamma(x), \log ||h(\gamma)(\dev(x))|| + t).\]
This preserves the affine structure and the radial vector field. 
The automorphism $\phi$ lifts to $\tilde \phi: \tilde \Sigma \ra \tilde \Sigma$ so that 
$\psi \circ \dev = \dev \circ \tilde \phi$ for $\psi \in \SL_{\pm}(n, \bR)$ where we define
\[\tilde \phi: \tilde \Sigma \times \bR \ra \tilde \Sigma \times \bR
\hbox{ by } 
\hat \phi(x, t) = (\tilde \phi(x),  \log ||\psi(\dev(x))|| + t).\] 
Then the result $\tilde \Sigma \times \bR/\langle \hat \phi, \pi_{1}(\Sigma)\rangle$   
 is homeomorphic to the mapping torus. 
We call this construction or the manifold the {\em generalized affine suspension}. 
If $\phi$ is of finite order, then the manifold is called a {\em Benz\'ecri suspension}. 

\begin{proposition} \label{prop:tcs} 
A connected compact radiant affine $n$-manifold 
 is a generalized affine suspension if and only if the radial flow has a total cross section.
\end{proposition}


\begin{corollary}[Barbot-Choi, Corollary A \cite{rdsv}] \label{cor:BC}
Let $M$ be a connected compact radiant affine $3$-manifold with empty or totally geodesic boundary, 
and the boundary is tangent to the radiant vector field. 
Then $M$ admits a total cross-section to the radiant flow. 
As a consequence, $M$ is affinely diffeomorphic to one of the following affine manifolds: 
\begin{itemize} 
\item a Benz\'ecri suspensions over a real projective surface of negative Euler characteristic with 
empty or geodesic boundary. 
\item a generalized affine suspension over a real projective sphere, a real projective plane, or 
a hemisphere, 
\item a generalized affine suspension over  a real projective torus {\rm (}Klein bottle{\rm )}, a real 
projective annulus {\rm (}M\"obius band{\rm )}  with geodesic boundary. 
\end{itemize} 
\end{corollary} 

There is a $6$-dimensional closed radiant affine manifold giving us a counter-example due to D. Fried \cite{Friednote} 

\begin{rem} 
We mention an error in \cite{rdsv} for Theorem A and Corollaries A and B. 
We state Corollary A in the corrected form above. 
We assume that not only that the holonomy group of the affine manifold $M$ fixes 
a common point but also we need that the boundary is tangent to the radial vector field. 
Proposition \ref{prop:rad} should fill in the gap since we just need to use the fact that 
radial lines foliate the universal cover. 
\end{rem}

\subsection{Affine $3$-manifolds with the infinite-cyclic holonomy groups.   } \label{sec:Hopf}

First, we will explore the affine Hopf manifolds. 

\begin{figure}[h]
	
	\begin{center}
		\includegraphics[height=6.5cm]{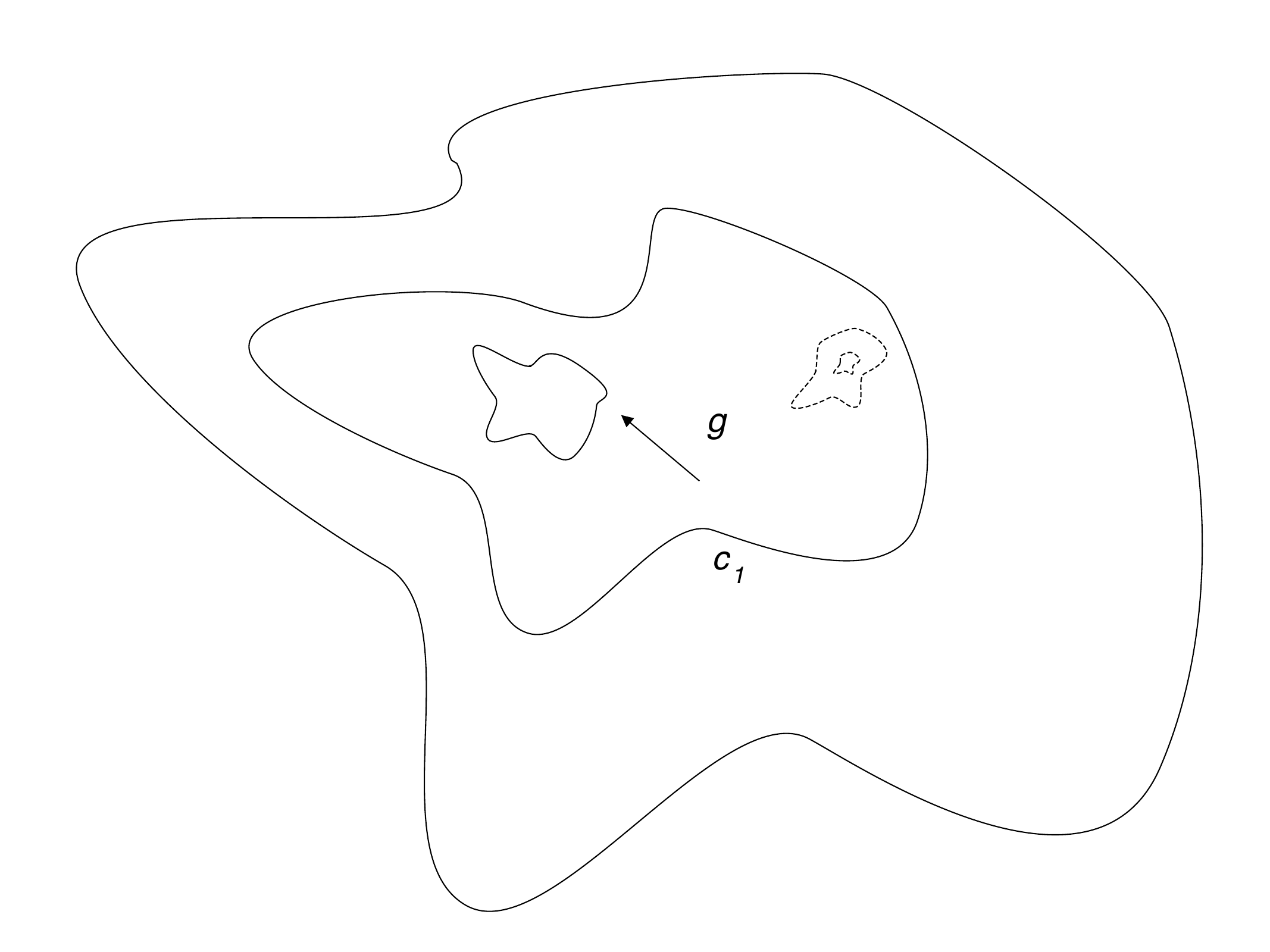}
	\end{center} 
	\caption{There must be an image of $c_1$ inside the component bounded by $c_1$.} 
	\label{fig:concave2}
\end{figure}

\begin{lemma} \label{lem:disj} 
Let $X$ be an open orientable manifold with a group $G$ acting on it properly discontinuously and cocompactly. 
Let $c_1$ be a compact connected submanifold where $X  - c_1$ has two open components, 
and let $U$ be a component. 
Then there exists infinitely many elements $g\in G$ so that $g(c_1) \subset U$. 
\end{lemma} 
\begin{proof} 
Let $x \in c_1$. 
Since the action of $G$ on $X$ is cocompact, 
there exists an infinite sequence $\{g_{i}\}$ of orientation-preserving 
$g_{i} \in G$ so that $g_{i}(x) \in U$. 
Since the action of $G$ on $X$ is properly discontinuous, $g_{i}(c_1) \cap c_1 = \emp$ except for 
finitely many $i$. We may choose $g_{i}$ so that $g_{i}(c_1)$ is a proper subset of $U$. 
\end{proof} 


\begin{proposition}\label{prop:Hopf} 
An affine Hopf $3$-manifold $M$ is homeomorphic to $\SI^{2}\times \SI^{1}$, $\bR P^{2}\times \SI^{1}$, 
or a nonorientable $\SI^{2}$-bundle over $\SI^{1}$. 
A half-Hopf $3$-manifold  $M$ is homeomorphic to a solid torus or a solid Klein bottle. 
\end{proposition}
\begin{proof} 
Let $M$ be an affine Hopf $3$-manifold. The universal cover is $\bR^{3}-\{O\}$ and hence 
$M$ does not contain any fake cells. 
We double cover it so that it has an infinite cyclic holonomy group
and call the double cover by $M'$. 
Let $g$ be the generator of the holonomy group. 
Each eigenvalue of a nonidentity element $g \in h(\pi_{1}(M))$ has either all norms $> 1$
or all norm $< 1$ by definition. 
By taking $g^{-1}$ if necessary, we assume that all the norms are $< 1$. 
Let $S$ be a unit sphere for a norm in Lemma \ref{lem:Katok}. 
By Lemma \ref{lem:Katok}, $S$ and $g(S)$ are disjoint.  
Then $S$ and $g(S)$ bound a compact space homeomorphic to 
$S \times I$. We introduce an equivalence relation $\sim$
where $x \sim y $ for $x\in S, y \in g(S)$ if $y = g(x)$. 
Thus, \[(\bR^{3} - \{O\})/\langle g \rangle \] is an $\SI^{2}$-bundle over $\SI^{1}$. 
Since $\Mod(\SI^2) = \bZ/2\bZ$ is a classical work of Smale \cite{Sm}, 
there exist only two homeomorphism types of $\SI^{2}$-bundle over $\SI^{1}$.

Now, $M$ is doubly or quadruply covered by $\SI^{2}\times \SI^{1}$. 
Since $-\Idd$ acts on $S$ above, and $\Mod(\bR P^{2}) = 1$,   
the proposition is proved. 

Suppose that $M$ is a half-Hopf $3$-manifold. Then we take a copy $M'$ of $M$ 
and glue $M$ with $M'$ at the boundary $\partial M$ and $\partial M'$ 
by a map induced by $-\Idd: \bR^{3}\ra \bR^{3}$ fixing $O$. 
Then the topology follows. 



\end{proof} 

See Section 3 of \cite{uaf} for the definition of the generalized affine suspension. 


\begin{theorem} \label{thm:cyclic}
Let $M$ be a connected compact affine $3$-manifold with empty or totally geodesic boundary  
and a virtually infinite-cyclic holonomy group 
whose infinite-order generator fixes a point in the affine space.
Also, the radial flow is tangent to the boundary. 
Then 
\begin{itemize}
\item $M$ is finitely covered by $\SI^2 \times \SI^1$ or $D^2 \times \SI^1$.
\item $M$ is a generalized affine suspension of a sphere, $\bR P^2$, or a \hyperlink{term-hem}{$2$-hemisphere}. 
\item If $M$ is closed, then $M$ is 
an affine Hopf $3$-manifold and is diffeomorphic to an $\SI^{2}$-bundle over $\SI^{1}$ or $\bR P^{2}\times \SI^{1}$. 
If $M$ has totally geodesic boundary, then $M$ is a half-Hopf manifold. 
\item Any $3$-manifold covered by an affine Hopf $3$-manifold or a half-Hopf $3$-manifold respectively is one also. 
\end{itemize} 
\end{theorem}
\begin{proof} 
We take a finite cover $N$ so that $N$ has an infinite cyclic fundamental group. 
By Theorem 5.2 of \cite{Hempel} and Lemma \ref{lem:3mfld}, $N$ has to be covered by $\SI^{2} \times \SI^{1}$ or $D^{2}\times \SI^{1}$ finitely. 
Therefore, the universal cover $\tilde M$ is neither complete affine nor bihedral. 

By taking a finite cover $N$ of $M$,  
we may assume that $h(\pi_1(N)) = \langle g \rangle$
and $g$ fixes a point $x$ in the affine space. 
Thus the holonomy group fixes a point $x$. 
Then $N$ is a radiant affine $3$-manifold by definition in \cite{rdsv}.  (See Section \ref{sub:ram}.)
Since the holonomy group is virtually infinite cyclic, 
the classification of such affine $3$-manifolds in 
Corollary \ref{cor:BC} (Corollary A in \cite{rdsv}) 
implies that $N$ is a generalized affine suspension of $\SI^{2}$, $\bR P^{2}$, or a $2$-hemisphere. 
To explain, $N$ admits a total cross-section by Theorem B of Barbot \cite{Barbot}. 
This means that $N$ and hence $M$ are covered by $\bR^{3}-\{x\}$ or $H-\{x\}$ for the closed half-space $H$ of $\bR^{3}$
for $x \in \partial H$. 

We now prove that when $M$ is closed, the only case is the affine Hopf $3$-manifold: 
$M$ is a generalized affine suspension of a real projective $2$-sphere or a real projective plane by the second item
of the conclusion of Corollary \ref{cor:BC}. 
In the first case,  $M$ has an infinite cyclic group as the deck transformation group
acting on $\bR^{3}-\{O\}$. By Proposition \ref{prop:hopf}, $M$ is an affine Hopf $3$-manifold. 
In the second case, the double cover of $M$ is an affine Hopf $3$-manifold. 
An order-two element $k$ centralizes the infinite cyclic group since the generator fixes a unique point in $\bR^{3}$. 
$\pi_{1}(M) $ is isomorphic to $\bZ \times \bZ_{2}$. 
Also, $k$ must act on a sphere in $\bR^{3} -\{O\}$ as an order two element, 
and, hence, $k = -\Idd$. Thus, $M$ is an affine Hopf $3$-manifold. 

When $M$ is a generalized affine suspension over a $2$-hemisphere, similar arguments apply to show that 
$M$ is a half-Hopf manifold. 


%
%

Any affine $3$-manifold covered by an affine Hopf $3$-manifold or a half-Hopf $3$-manifold satisfies the premises of the theorem. 
Thus, it is an affine Hopf $3$-manifold or a half-Hopf $3$-manifold. 


\end{proof}

\begin{corollary} \label{cor:cyclic} 
Let $M$ be a connected closed real projective $3$-manifold. 
Suppose that $M_h$ is a domain $\Omega$   
in $\SI^3$ containing a smoothly embedded sphere $S^2$ of codimension one 
bounding an $3$-ball $B^3$ in $\bR^3$ for an affine subspace $\bR^3 \subset \SI^3$, and 
$B^3$ is not contained in $\Omega$. 
Then $M$ is projectively diffeomorphic to an affine Hopf $3$-manifold. 
\end{corollary} 
\begin{proof} 
A component of $\SI^3 - S^2$ is the $3$-cell $B^3$ in an affine patch. So, $B^3$ is in a properly convex domain. 
Since $M_{h}$ covers a compact manifold, there exists a deck transformation $g$ so that 
$g(S^2) \subset B^3\cap \Omega$ by Lemma \ref{lem:disj}. 
Then $g(B^3) \subset B^3$ since  the outside component of $\SI^2 - g(S^2)$ is not contained in a properly convex domain. 
By the Brouwer fixed-point theorem, $g$ fixes a point in the interior of 
$B^3$. 
Since $g$ sends the closure of the neighborhood of $x$ bounded by $S^2$ 
into itself, a simple argument shows that the one-dimensional subspace in $\bR^{4}$ 
in the direction of $x$ has the largest norm of multiplicity one for 
the linear representative of $g$. Thus, 
$g$ acts on a hyperspace $W$ of $\SI^3$ disjoint from $x$
and acts as an affine transformation on  the affine space $\bH^{o}$ bounded by $W$
and containing $x$. 
Here, $B^3 \subset \bH^o$ clearly. 
Proposition \ref{prop:diskf} implies 
that $x$ is an attracting fixed point of $g$ on $\bH^{o}$. 



By Lemma \ref{lem:residual}, 
there is a cover $M'$ of $M$ by taking a finite-index normal subgroup $\Gamma_{h}'$ of $\Gamma_{h}$
so that $p|S^2$ is an embedding to a sphere $\hat S$
for the covering map $p:\Omega \ra M'$. 
Furthermore, we may assume that $\Gamma_{h}'$ is orientation preserving. 

Let $g^{i_{0}}$ be the power of $g$ in $\Gamma'_{h}$ with least $i_{0}, i_{0}>0$.
Now, $g$ is a contracting 
linear map with a fixed point $x \in B^3$, and 
\[g(S^2) \subset B^3 \hbox{ and } B^3 \subset x\ast S^2.\] 
Thus, for sufficiently large $i_{0}$,  we obtain 
\[g^{i_{0}}(S^2) \subset B^3 \hbox{ and }  g^{i_{0}}(B^3 \cap \Omega) \subset B^3\cap \Omega.\] 
 as in the beginning of the proof. 

For any $k \in \Gamma'_{h}$ so that $k(S^2) \subset B^3$, 
we have $k(B^3) \subset B^3$ as in the beginning of the proof. 
Thus any generic open arc $\alpha$ connecting $x \in S^2$ to $g^{i_{0}}(x)$ in $B^3\cap \Omega - g^{i_{0}}(\clo(B^3))$
meets copies of $B^2$ under $\Gamma'_{h}$ other than $g^{i_{0}}(B^3)$ 
\begin{itemize} 
\item in a compact interval disjoint from $x$ and $g^{i_{0}}(x)$ or
\item in the interval containing $g^{i_{0}}(x)$ but not $x$. 
\end{itemize}
Lets call $f$ the number of times the second case happens. 
We may assume that the intersection number at $x$ with $\alpha$ is $+1$. 
Then the oriented-intersection number of the image of $\alpha$ in $M'$ meeting
$\hat S$ is $f+1 > 0$. 
Thus, $\hat S$ is a nonseparating sphere. 



Since $\hat S$ is nonseparating, we take a transversal embedded arc $I'$ in $M'$ connecting a point of $\hat S$ to itself and 
disjoint from $\hat S$ in the interior. We take an $\eps$-neighborhood $N$ of $I'\cup \hat S$ and 
let $S''$ denote the boundary sphere of the neighborhood.  (See Lemma 3.8 of \cite{Hempel} for the construction.) 
Since $S''$ bounds a neighborhood, it is a separating $2$-sphere. 
We choose sufficiently small $\eps$ so that $S''$ is homeomorphic to a sphere. 
Let $I$ be the lift of $I'$ starting from $S^2$. 
Let $k$ be the deck transformation so that $y$ and $k(y)$ are endpoints of $I$. 
Then $S''$ lifts to a sphere in $\Omega$ that is a boundary component $S'$ of the inverse image $N''$ of $N$, which is 
a neighborhood of $\bigcup_{n\in \bZ} k^3(S^2 \cup I)$. 

Again, $k$ satisfies the properties of $g$ above. 
By change of notation, let $x$ denote to the attracting fixed point of $k$ of an affine space to be denote by $\bH^{o}$. 

Since a sequence $\{k^{j}(S^2)\}$ of compact sets 
geometrically converges to $\partial \bH$ as $j \ra -\infty$,  
and $S^2$ is compact, $\Omega$ is disjoint from $\partial \bH$ by
the properness of the action of $\langle k \rangle$. 
Hence, $\Omega \subset \bH^o$. 
Since $x$ is a fixed point of $k$, we obtain $\Omega \subset \bH^{o} -\{x\}$.

Since $\Omega \subset \bR^3$, the sphere $S'$ is a subset of $\bR^3$. 
By Theorem 1.2 of Wu \cite{Wu}, $S''$ bounds a $3$-ball $B''$ in $M'$, 
and a $3$-ball $B'$ bounded by $S'$ in $\Omega=M_{h}$ 
embeds onto $B''$ under $p$. 
Since  $\partial B'$ is a component of  $\partial N''$,
we cannot have $B' \subset N''$ while $N''$ is not compact. 
Hence, the interior of $B'$ is disjoint from that of $N''$. 

Taking a union of $B'$ with $N''$, 
we obtain that the domain $D$ in $\bH^{o}$ bounded by $S^2$ and $k(S^2)$ is in $\Omega$. 
Since 
\[D\subset B' \cup N \subset\Omega,\]
we obtain $\Omega \supset \bigcup_{n\in \bZ} k^3(D)$. 

By the generalized Schoenflies theorem 
$k^{-n}(S^2)$ and $k^3(S^2)$ bound a region in $\bH^{o}$ homeomorphic to $S^2\times I$, and 
\[k^{-n}(S^2) \ra \partial \bH \hbox{ and } k^3(S^2) \ra \{x\} \hbox{ as } n \ra \infty\]
in the Hausdorff convergence sense.  
It follows that $\bigcup_{n\in \bZ} k^3(D) = \bH^{o} -\{x\}$. 
Hence, $\Omega = \bH^{o} - \{x\}$. 


Now, Theorem \ref{thm:cyclic} shows that $\Omega/\langle g^{{i_{0}}} \rangle$ is 
an affine Hopf $3$-manifold, a compact manifold. 
Therefore, $M$ is finitely covered by an affine Hopf-$3$-manifold. 
Theorem \ref{thm:cyclic} implies the result. 

\end{proof}

\subsection{Convex concave decomposition of real projective $3$-manifolds} \label{subsec:convconc}

\subsubsection{Kuiper completions} \label{subsub:Kuiper} 

The immersion $\dev_h$ induces a Riemannian $\mu$-metric on $M_h$ from the standard Riemannian metric $\mu$ on $\SI^3$. 
This gives us a path-metric to be denoted by $\bdd$ on $M_h$. (More precisely $\bdd_h$ but we omit $J$ here.)  
Recall from \cite{psconv}
the Cauchy completion $\che M_h$ of $M_h$ with this path-metric 
is called the \hypertarget{term-kcp}{{\em Kuiper completion}} of $M_h$. (This metric is quasi-isometrically defined by $\dev_h$
and hence the topology is independent of the choice of $\dev_h$.)
The \hypertarget{term-ids}{{\em ideal set}} is $M_{h, \infty} := \che M_h - M_h$, which is in general not empty. 
The immersion $\dev_h$ extends to a continuous map. We use $\dev_h$ as the notation for the extended map as well. 
If $M$ is an affine $3$-manifold, we define $M$ as a real projective $3$-manifold 
and $\che M_h$ as above. $\Gamma_h$ acts on $M_h$ and $M_{h, \infty}$ possibly with fixed points in $M_{h, \infty}$. 

For a compact convex subset $K$ of $\che M$ so that $\dev_h|K$ is an embedding, we define 
$\partial K$ to be the subset corresponding to $\partial \dev_h(K)$. 
If $\dev_h(K)$ is a compact domain in a subspace of $\SI^3$, 
then we define $K^o$ as the subset of $K$ that is the inverse image of the manifold interior of $\dev_h(K)$.
An \hypertarget{term-hem}{{\em $i$-hemisphere}} in $\che M_h$ 
is a compact subset $H$ so that $\dev_h|H$ is an embedding to an $i$-hemisphere, $1 \leq i \leq 3$.
For $3$-hemisphere, we require $H^{o}\subset M_{h}$. 
A \hypertarget{term-bih}{{\em $3$-bihedron}} 
in $\che M_h$ is a compact subset $B$ so that $B^{o}\subset M_{h}$ and $\dev_h|B$ is an embedding to a compact convex set $K$ so that 
$\partial K$ is the union of two $2$-hemispheres with the identical boundary great circle. 

\subsubsection{$2$-convexity and covers} 
A {\em tetrahedron} or $3$-simplex is a convex hull of four points in general position in an affine space $\bR^3$.
A real projective $3$-manifold $M$ is {\em $2$-convex} if every projective map 
$f:T^o \cup F_1 \cup F_2 \cup F_3 \ra M$ for a tetrahedron $T$ with faces $F_i, i=0, 1, 2, 3$, 
extends to one from $T$. 


A {\em tetrahedron} in $\che M_h$ is a compact subset $T$ 
so that $\dev_h|T$ is an embedding to a tetrahedron in an affine space in $\SI^3$.
A {\em face} of $T$ is a corresponding subset of $\dev_h(T)$.

\begin{proposition}[ Proposition 4.2 of \cite{psconv}]\label{prop:tet} 
$M$ is $2$-convex if and only if for every tetrahedron $T$ in $\che M_h$ with faces $F_i$, $i=0, 1, 2, 3$, such that
$T^o \cup F_1 \cup F_2 \cup F_3 \subset M_h$, $T$ is a subset of $M_h$.
\end{proposition}


\subsubsection{Crescents and two-faced submanifolds} \label{subsec:crcaf}

A \hypertarget{term-hcr}{{\em hemispherical $3$-crescent}} is a $3$-hemisphere $H$ in $\che M_h$ so that 
$H^{o}\subset M_{h}$ and 
the $2$-hemisphere in $\partial H$ is 
a subset of the ideal set $M_{h, \infty}$. We define $\alpha_R$ for a hemispherical $3$-crescent $R$ to be the union of 
all open $2$-hemispheres in $\partial R \cap M_{h, \infty}$. We define $I_R = \partial R - \alpha_R$.

By Proposition 6.2 of \cite{psconv} or by its $\che M_h$-version, 
given two hemispherical $3$-crescents $R$ and $S$ in $\che M_h$, the exactly one of following hold: 
\begin{itemize}
\item $R \cap S \cap M_{h}= \emp$,
\item  $R = S$, or 
\item $R \cap S \cap M_{h}$ is a union of common components of $I_R \cap M_h$ and 
$I_S \cap M_h$. 
\end{itemize} 
The components of $I_R \cap M_h$ as in the last case are called {\em copied 
components} of $I_R \cap M_h$. 
The union of all copied components in $M_h$, a {\em pre-two-faced submanifold of type I}, is totally geodesic and 
covers a compact embedded totally geodesic $2$-dimensional submanifold in $M^o_h$
by Proposition 6.4 of \cite{psconv}. 
The submanifold is called the \hypertarget{term-tsfm}{{\em two-faced submanifold of type I}} (arising from hemispherical $3$-crescents).  
(It is possible that the needed results of \cite{psconv} are true when the manifold-boundary $\partial M$ is convex.
This is not proved there.) Note it is possible that the two-faced submanifold of type I may be empty, i.e., does not exist at all. 

The \hypertarget{term-spl}{{\em splitting along}} a submanifold $A$ is given by the Cauchy completion $M^s$ of $M-A$ 
of the path metric obtained by using an ordinary Riemannian metric on $M$ and restricting to $M-A$. 



\begin{hyp}\label{hyp:nohem}
We assume that  a bihedral $3$-crescent in $\che M_{h}$ is defined if 
there is no hemispherical $3$-crescent. 
\end{hyp} 
In other words, we shall talk about bihedral $3$-crescent when $\che M_{h}$ has no hemispherical $3$-crescent.

A \hypertarget{term-bcr}{{\em bihedral $3$-crescent}} is a $3$-bihedron $B$ in $\che M_h$ so that 
$B^{o}\subset M_{h}$ and 
a $2$-hemisphere in $\partial B$ is a subset of $M_{h, \infty}$. 
(We require that these are not contained in a hemispherical $3$-crescent.)
For a bihedral $3$-crescent $R$,  we define $\alpha_R$ as the open $2$-hemisphere in $\partial R \cap M_{h, \infty}$. 
We define $I_R := \partial R - \alpha_R$, a $2$-hemisphere. 
For a $3$-crescent $R$, we define the interior of $R$ as $R^o = R - I_R - \alpha_R$.

We say that two $3$-crescents $R$ and $S$ \hypertarget{term-ovl}{{\em overlap}} if $R^o \cap S \ne \emp$, or equivalently 
$R^o \cap S^o \ne \emp$. 
We say that $R \sim S$ if 
there exists a sequence of $3$-crescents $R_1 = R, R_2, \cdots, R_n = S$ 
where $R_i \cap R_{i+1}^o \ne \emp$ for $i=1, \dots, n-1$. 

We say that two bihedral $3$-crescents $R$ and $S$ intersect {\em transversally} if 
\begin{itemize}
\item $I_{S}\cap I_{R}$ is a segment with endpoints in $\partial I_{S}$ and $\partial I_{R}$, 
\item $I_{S}\cap R$ is the closure of a component of $I_{S} - I_{R}$, and 
\item $R \cap S$ is the closure of a component of $R - I_{S}$.
\end{itemize} 
In this case, $\alpha_{S} \cup \alpha_{R}$ is a union of two open $2$-hemispheres meeting 
at an open convex disk $\alpha_{S}\cap \alpha_{R}$. 
Thus, they {\em extend} each other. 
(See Chapter 5 of \cite{psconv}.) 

\begin{proposition} \label{prop:transversal} 
We assume as in Theorem \ref{thm:main1}. 
Suppose that two bihedral $3$-crescents $R$ and $S$ in $\che M_{h}$ overlap. Then $R$ and $S$ either intersect transversally or 
$R \subset S$ or $S \subset R$. Moreover, $\dev_{h}| R \cup S$ is a homeomorphism to its image 
$\dev_{h}(R) \cup \dev_{h}(S)$ where $\dev_{h}(\alpha_{R})$ and $\dev_{h}(\alpha_{S})$ are $2$-hemispheres in 
the boundary of a $3$-hemisphere $H$. 
\end{proposition} 
\begin{proof} 
This is a restatement of Theorem 5.4 and Corollary 5.8 of \cite{psconv}.
\end{proof} 


From now on, assume that there is no hemispherical $3$-crescent in $\che M_{h}$.
We define as in Chapter 7 of \cite{psconv} 
\begin{align} \label{eqn:defL}
 \Lambda(R) &:= \bigcup_{S \sim R} S, \quad \quad \quad \delta_\infty \Lambda(R) := \bigcup_{S \sim R} \alpha_S, \nonumber \\ 
  \Lambda_1(R) &:= \bigcup_{S\sim R} (S- I_S), \quad \delta_\infty \Lambda_1(R) :=  \delta_\infty \Lambda(R). 
 \end{align}
We showed in Chapter 7 of \cite{psconv}
 $\dev_h|\Lambda(R)$ maps into a $3$-hemisphere $H$ and 
 $\dev_h| \delta_\infty \Lambda(R)$ is an injective local homeomorphism to $\partial H$  (see also Corollary 5.8 of \cite{psconv}). 
 

Given a subset $A$ of $\che M_h$, we define $\hbox{int} A$ to be the interior of $A$ in $\che M_h$. 
We define $\hbox{bd} A$ to be the topological boundary of $A$ in $\che M_h$. 
By Lemma 7.4 of \cite{psconv}, there are three possibilities: 
\begin{itemize} 
\item if $\inte \Lambda(R)  \cap \Lambda(S) \cap M_h \ne \emp$ for two bihedral $3$-crescents $R$ and $S$, 
then $\Lambda(R) = \Lambda(S)$, 
\item $\Lambda(R) \cap \Lambda(S) \cap M_{h} = \emp$, or
\item $\Lambda(R) \cap \Lambda(S) \cap M_h \subset \Bd \Lambda(R) \cap M_h \cap \Bd \Lambda(S) \cap M_h$. 
\end{itemize} 
In the third case, the intersection is a union of common components of $\Bd \Lambda(R) \cap M_h$ and $\Bd \Lambda(S) \cap M_h$. 
We call such components {\em copied components}. These are totally geodesic and properly embedded in $M_{h}$. 
The union of all copied components in $M_h$, a {\em pre-two-faced submanifold of type II}, 
covers a compact embedded totally geodesic $2$-dimensional submanifold in 
$M^o$ by Proposition 7.6 of \cite{psconv}. 
The submanifold is called the \hypertarget{term-tfsm}{{\em two-faced submanifold of type II}} (arising from bihedral $3$-crescents). 

%


\subsubsection{Concave affine manifolds after splitting} \label{subsub:split} 

Let $M^s$ denote the $3$-manifold obtained from $M$ by splitting along the two-faced submanifold of type I. 
A cover $M^s_{h}$ of $M^s$ can be obtained 
by splitting along the preimage of the two-faced submanifold of type I
 in $M_h$ and taking a component for every component of $M^s$ and taking the union of these. 
For each component $A$ of $M^s_h$, let $\Gamma_A$ denote the subgroup of $\Gamma_h$ acting on 
$A^o$. Then $\Gamma_A$ extends to a deck transformation group of $A$. 
We define $\Gamma^s_h$ the product group 
\[\prod_{A \in \mathcal{C}} \Gamma_A  \hbox{ for the set } \mathcal{C} \hbox{ of chosen components in } M^s_h.\] 
 Again $M^s_h$ has a developing map $\dev^s_h: M^s_h \ra \SI^3$, an immersion, and 
 $M^s_h \ra M^s$ is a holonomy cover with the deck transformation group $\Gamma^s_h$. 
 There is a map $M^s_{ h} \ra \che M$ by identifying along the splitting submanifolds. 
We can easily see that the \hyperlink{term-kcp}{Kuiper completion} $\che M^s_{ h}$ contains the hemispherical $3$-crescent if and only if $\che M_h$ does. 
Also, the set of hemispherical $3$-crescents of $\che M^s_{ h}$ is mapped in a one-to-one manner to the set of those in $\che M_{h}$ 
by taking the interior of the hemispherical $3$-crescent in $\che M^{s}_{h}$ and sending it to $\che M_{h}$ and taking the closure. 
Now $\che M^s_{h}$ does not have any copied components. 
(See Chapter 8 of \cite{psconv}).

\begin{definition} \label{defn:cafI} 
A connected compact real projective manifold with totally geodesic boundary covered by $R \cap M^{s}_{h}$
a hemispherical $3$-crescent $R$ is said to be a \hypertarget{term-camI}{{\em concave affine manifold of type I}} in $M^{s}$. 
\end{definition} 



Let $\mathcal{H}$ be the set of all hemispherical $3$-crescents in $M^{s}_{h}$. 
The union $\bigcup_{R \in \mathcal{H}}R\cap M^{s}_{h}$
covers a finite union $K$ of mutually disjoint concave affine manifolds of type I in $M^{s}$. 
Then $M^{s} - K^{o} = M^{(1)}$ is a compact real projective manifold with convex boundary. 
The cover $M^{(1)}_{h} $ of $M^{(1)}$ is $M^{s}_{h}$ with all points of hemispherical $3$-crescents removed from it. 
Then $\che M^{(1)}_{h}$ has no hemispherical $3$-crescent. (See p. 80--81 of \cite{psconv}.)

Now, we look at $M^{(1)}$ only. 
We split $M^{(1)}$ along the two-faced submanifold of type II if it exists. Let $M^{(1)s}$ denote the result of the splitting. 
Also, the set of bihedral $3$-crescents of $\che M^{(1)}$ is mapped in a one-to-one manner to the set of those in $\che M^{(1)s}$ 
by taking the interior of the bihedral $3$-crescent and sending it to $M^{(1)s}$ and taking the closure. (See Chapter 8 of \cite{psconv}).
Now $\che M^{(1)s}_{h}$ does not have any copied components. For a bihedral $3$-crescent $R$ in $\che M^{(1)s}_h$, 
$\Lambda(R) \cap M^{(1)s}_{h}$ covers a compact $3$-manifold with concave boundary in $M^{(1)s}_{h}$. 
(See p. 81--82 of \cite{psconv}.)

\begin{definition}\label{defn:cafII} 
Suppose that $\che M_{h}$ does not contain a hemispherical $3$-crescent. (See Hypothesis \ref{hyp:nohem}.)
Let $R$ be a bihedral $3$-crescent in $\che M_{h}$. 
If $\Lambda(R)\cap M_{h}$ covers a compact real projective submanifold $N$,
then $N$ is called a \hypertarget{term-camII}{{\em concave affine manifold of type II}}. 
\end{definition} 
See Chapter 8 of \cite{psconv} as a reference of results stated here.

\section{Concave affine $3$-manifolds} \label{sec:concave}

In this section, we will prove Theorems \ref{thm:tfaced}, \ref{thm:caffI} and \ref{thm:caffII}. 
The first one shows that the non-$\pi_{1}$-injective two-faced submanifolds cannot happen in general. 
In the second and third ones, we showed that 
a concave affine manifold with compressible boundary contains a \hyperlink{term-tsp}{toral $\pi$-submanifold}. 

Given an embedded surface $\Sigma$ in a $3$-manifold $M$ that is either on the boundary of $M$ or is two-sided, 
$\Sigma$ is {\em incompressible} into $M$ if 
$\pi_1(\Sigma)$ injects into $\pi_1(M)$. Otherwise, $\Sigma$ is said to be {\em compressible}.
A simple closed curve in $\Sigma$ is {\em essential} if it is not null-homotopic in $\Sigma$. 
A compressible surface always has an essential simple closed curve that is the boundary of an embedded disk by
Dehn's Lemma. 

\subsection{A concave affine manifold has no sphere boundary component}




\begin{theorem}\label{thm:bdsph}
Let $N$ be a concave affine manifold of type I or II. 
Then no component of $\partial N$ is covered by a sphere.
\end{theorem} 
\begin{proof} 
If $N$ is a concave affine manifold of type I, then $N$ is covered by $\tilde N= R^{o} \cup I_{R}\cap N_{h}$
for a hemispherical $3$-crescent $R$.  Since $I_{R} \cap N_{h}$ is an open surface in $I_{R}$, 
the conclusion follows. 

Suppose now that there is no hemispherical $3$-crescent in $M_{h}$. (See Hypothesis \ref{hyp:nohem}.)
Then $N$ is covered by $\Lambda(R) \cap M_{h}$ for a bihedral $3$-crescent $R$. 
Suppose that a component $A$ of $\Bd \Lambda(R) \cap M_h$ is a sphere. 
We know that $A$ is mapped into a convex surface in $M-N^o$ under the covering map. 
If $A$ is totally geodesic, then $A$ is tangent to $I_{S} \cap M_{h}$ for a crescent $S$ in $\Lambda(R)$. 
Hence, $A$ is a subset of $I_S \cap M_h$, each component of which is not compact. 
This is a contradiction.

Suppose that each point $x$ of $A$ has some open geodesic segment in $A$ containing $x$. 
Since $A$ is convex, $x$ is on a unique maximal geodesic in $A$ or is in a $2$-dimensional totally geodesic surface in $A$.
Since $A$ is convex, a geodesic segment in $A$ must end at the boundary of $A$. 
This implies that $A$ is not compact, a contradiction. 

Hence, there must be a point $y$ where $A$ is strictly concave. 
This contradicts Theorem \ref{thm:nscc}. 

\end{proof}

\subsection{Non-$\pi_{1}$-injective two-faced submanifolds}\label{sub:tfaced}

\begin{lemma} \label{lem:simplyc} 
Let $\tilde A_{1}$ be a properly embedded surface in $M_{h}$ covering a compact surface $A_{1}$. 
If $\tilde A_{1}$ is a disk, then $A_{1}\ra M$ is $\pi_{1}$-injective. 
\end{lemma} 
\begin{proof} 
The deck transformation group $\Gamma_{\tilde A_{1}}$ acting on $\tilde A_{1}$
injects into the deck transformation group $\Gamma_{h}$. 
\end{proof}

\begin{theorem} \label{thm:tfaced} 
Suppose that a connected compact real projective $3$-manifold $M$ with empty or convex boundary
and $M$ is neither complete affine nor bihedral. 
Suppose that $M$ contains the \hyperlink{term-tfsm}{two-faced submanifold} $S$ in $M$. 
Then either $S$ is $\pi_{1}$-injective in $M$ or $M$ is an affine Hopf $3$-manifold. 
\end{theorem}

This implies that two-faced submanifolds are $\pi_{1}$-injective unless $M$ is an affine Hopf $3$-manifold.



(I) Let $A_1$ denote a component of the two-faced submanifold of type I. 
Suppose that $A_1 $ is covered by a component $\tilde A_1$ of $I_R \cap M_h$ for a hemispherical $3$-crescent $R$. 
If $\tilde A_1$ is simply connected and planar, then $\tilde A_{1}$ is a disk. 
By Lemma \ref{lem:simplyc}, 
$A_1$ is $\pi_{1}$-injective in $M$, and we are done here. 

Let $\Gamma_1$ denote the deck transformation group of $\tilde A_1$ in $\Gamma_h$ so that $\tilde A_1/\Gamma_1$ is compact
and diffeomorphic to $A_1$. 
Now assume that $A_{1}$ is non-$\pi_{1}$-injective in $M$. 
By the above paragraph, the planar surface
$\tilde A_1$ contains a simple closed curve $c_1$ not bounding a disk in $\tilde A_1$. 

By Corollary \ref{cor:cyclic}, $\tilde A_1$ is projectively 
diffeomorphic to $\bR^2 -\{O\}$. Since $\tilde A_1 \subset I_R^o$, 
we  obtain \[\tilde A_1 = I_R - \{x\} \hbox{ for } x\in I_{R}^{o}.\] 
Since $\tilde A_1$ covers a component of the two-faced submanifold, 
$\tilde A_1$ is a component of $I_S \cap M_h$ for a hemispherical $3$-crescent $S$ where 
\[R \cap S \cap M_h = \tilde A_1.\] 

Since $\clo(\alpha_S) \cup \clo(\alpha_R) \subset M_\infty$ bounds the compact domain $R \cup S$ in $\che M$, 
we obtain $R^o \cup S^o \cup \tilde A_1 = M_h$. 
Now, $\dev_{h}| R^{o} \cup I_R^o -\{x\} $ and $\dev_{h}| S^{o} \cup I_S^o -\{x\})$ are homeomorphisms to their images. 
Thus, $\dev_{h}| M_{h}$ is a homeomorphism to the image 
\[\dev_{h}(R)^{o}\cup \dev_{h}(S)^{o} \cup \dev_{h}(I_{R}^{o}) - \dev_{h}(x).\]  
Since $\dev_{h}(x)$ is an isolated boundary point, 
by Corollary \ref{cor:cyclic}, we are finished in this case.  


(II) Let $A_1$ denote a component of a two-faced submanifold of type II in $M$ that is non-$\pi_{1}$-injective. 
Now, we assume that $\che M_{h}$ has no hemispherical $3$-crescent. (See Hypothesis \ref{hyp:nohem}.)
Then as in case (I), its cover $\tilde A_1$ is a component of $I_R \cap M_{h}$ 
containing a simple closed curve not contractible in 
$\tilde A_1$ for a bihedral $3$-crescent $R$. 
By Corollary \ref{cor:cyclic}, 
we obtain that $\tilde A_1 = I_R^o -\{x\}$ for a bihedral $3$-crescent $R$. 

Since $\tilde A_1$ is in a pre-two-sided submanifold, we obtain that $I_R \subset I_S$ 
for another bihedral $3$-crescent $S$ so that $S^o \cap R^o =\emp$. 
It follows that 
\[I_R^o -\{x\} = I_S^o -\{x\} \hbox{ and hence } I_R = I_S.\] 
Since $\clo(\alpha_R) \cup \clo(\alpha_S) \subset M_{h, \infty}$ forms the boundary of $R \cup S$, 
and $M_h$ is disjoint from it, 
\[M_h = R^o \cup S^o \cup I_R^o -\{x\}\] holds. 
Hence, $\dev_{h}$ is an embedding to 
$\dev_{h}(R^{o}) \cup \dev_{h}(S^{o}) \cup \dev_{h}(I_R^o -\{x\} )$. 
Since $\dev_{h}(x)$ is an isolated boundary point, 
Corollary \ref{cor:cyclic} implies the result in this case.




\begin{lemma}\label{lem:nobd} 
Let $\Omega_1$ be an open surface in $M_h$ with $\dev_h(\Omega_1)$ 
bounded in an affine space $H^o$ for a \hyperlink{term-shem}{$2$- or $3$-hemisphere} $H$. 
We assume that $H$ is the minimal dimensional hemisphere containing $S$. 
Suppose that a discrete group $G \subset \Gamma_h$ acts properly discontinuously and freely on $\Omega_1$, and $h_h|G$ is injective.
Moreover, $G$ acts on $H$. 
Then $\Omega_1/G$ is noncompact. 
\end{lemma} 
\begin{proof} 
Suppose that there exists $G$ so that $\Omega_1/G$ is compact. 
Since $G$ acts on $H$, $G$ acts as a group of affine transformations on the affine $2$- or $3$-space $H^o$.
Let $F$ be the compact fundamental domain of $\Omega_{1}$. 
The closure $\clo(\dev_h(\Omega_1))$ is a compact bounded subset of $H^{o}$. 
The convex hull $C_1$ of $\clo(\dev_h(\Omega_1))$ is a bounded subset of $H^{o}$ also,
and $G$ acts on it and its center of mass $m$, and hence $h_{h}(G)$ is a group of bounded  affine transformations
fixing $m$. 
We choose a $h_{h}(G)$-invariant Euclidean metric $d_{H^{o}}$ on $H^{o}$. 
Let $U$ be an open $\eps$-$d_{H^{o}}$-neighborhood of $F$ in $S$. 
We choose sufficiently small $\epsilon$ so that $U \subset \Omega$. 

Since $\Omega$ is open, there exists a sequence $\{y_i\}$ exiting all compact sets in $\Omega$ eventually.
There exists $g_{i} \in G$ such that $g_{i}(y_{i}) \in F$. 
By taking a subsequence, we may assume $\dev_{h}(y_{i}) \ra y \in S$ and $y $ is in the boundary of $\dev_{h}(\Omega)$, 
i.e., $y \not\in \dev_{h}(\Omega)$. 
Then $g_{i}^{-1}(F) \ni y_{i}$. 
Since $\dev_{h}(y_{i}) \ra y$, $h_{h}(g_{i}^{-1})$ is an isometry group fixing $m$, and $S$, $S \ni y$, is properly embedded, 
it follows that 
\[\dev_{h}(\Omega) \supset \dev_{h}(g_{i}^{-1}(U^{o})) = h_{h}(g_{i}^{-1})(\dev_{h}(U^{o})) \ni y \hbox{ for sufficiently large } i,\]
which is a contradiction. 
\end{proof}

\subsection{Concave affine manifolds and toral $\pi$-manifolds}\label{subsec:caff}

\begin{definition}\label{defn:opposite}
Suppose that $\che M_{h}$ contains two crescents $S_1$ and $S_2$ so that 
$I_{S_1} \cap M_h$ and $I_{S_2} \cap M_h$ 
intersect and are tangent but $\dev_h(S_1)^o \cap \dev_h(S_2)^o =\emp$. 
In this case $S_1$ and $S_2$ are said to be \hypertarget{term-opp}{{\em opposite}}.
\end{definition}

\begin{definition}\label{defn:toralpi} 
Suppose that a compact connected real projective manifold $M$ is neither complete affine nor bihedral,
and let $M_{h}$ be the holonomy cover of $M$. 
Assume that $M$ has no two-faced submanifolds. 
Let $R$ be a hemispheric $3$-crescent with $I_R\cap M_h = I_R^o -\{x\}$ for $x \in I_R^o$.
Then a compact submanifold $P$ covered by $R \cap M_h$ is called a {\em toral $\pi$-submanifold of type I}. 

Suppose that $\che M_{h}$ has no hemispheric crescent. (See Hypothesis \ref{hyp:nohem}.)
Given $\Lambda(R)$ for a bihedral $3$-crescent $R$, 
we define the set $C_{R, x}$ as follows: 
Suppose that for some $x \in \SI^{3}$, we define
\[ C_{R, x} := \{ R' | R' \sim R, \exists g \in \Gamma_{h}, g(R) = R, h_{h}(g)(x) = x,  \dev_{h}(I_{R'}^{o}) \ni x \} \ne \emp. \]
Let $\Lambda'(R)$ be $\bigcup_{R'\in C_{R, x}} R'$ whenever $C_{R, x}$ is not empty and 
\[\delta_\infty \Lambda'(R) := \bigcup_{S\in C_{R, x}} \alpha_{S}.\] 
Then $\Lambda'(R)$ develops into a $3$-hemisphere $H$, and $\delta_\infty \Lambda'(R)$ 
develops to an open disk in $\partial H$ for a $3$-hemisphere $H$ by Chapter 7 of \cite{psconv}. 
Suppose that $\Lambda'(R)\cap M_h$ covers a compact radiant affine $3$-manifold $P$ with compressible boundary. 
Then $P$ is said to be a \hypertarget{term-tsp}{{\em toral $\pi$-submanifold of type II}}.
\end{definition}


Theorems \ref{thm:caffI} and \ref{thm:caffII} characterize the 
concave affine $3$-manifolds with compressible boundary. One consequence is that 
the fundamental group is virtually infinite cyclic. 

\begin{theorem} \label{thm:caffI}
Let $N$ be a concave affine $3$-manifold with nonempty boundary $\partial N$ 
in a connected compact real projective $3$-manifold $M$ with empty or convex boundary.
Suppose that $M$ is neither complete affine nor bihedral. 
Assume that $M$ has no \hyperlink{term-tfsm}{two-faced submanifold of type I}. 
Let $M_{h}$ be the holonomy cover of $M$.
Suppose that $N$ is a concave affine $3$-manifold of type I with compressible boundary $\partial N$. 
Then one of the following holds\,{\rm :} 
\begin{itemize}
\item $M$ is an affine Hopf $3$-manifold,  or
\item $N$ has a unique boundary component $A$ compressible into $N$, and 
$N$ is a toral $\pi$-submanifold $P$ of type I. 
\end{itemize} 
\end{theorem} 

\begin{theorem} \label{thm:caffII}
Let $N$ be a concave affine $3$-manifold with nonempty boundary $\partial N$ 
in a connected compact real projective $3$-manifold $M$ with empty or convex boundary.
Suppose that $M$ is neither complete affine nor bihedral. 
Let $M_{h}$ be the holonomy cover of $M$. 
Suppose that $\che M_{h}$ has no hemispherical $3$-crescent. {\rm (}See Hypothesis \ref{hyp:nohem}.{\rm )}
Assume that $M$ has no two-faced submanifold of type II. 
Suppose that $N$ is a concave affine $3$-manifold of type II with compressible boundary $\partial N$. 
Then one of the following holds\,{\rm :} 
\begin{itemize}
\item $M$ is an affine Hopf $3$-manifold, or
\item $N$ has a unique boundary component $A$ compressible into $N$, and 
$N$ contains a maximal toral $\pi$-submanifold $P$ of type II. 
Furthermore, the following holds\,{\rm :}
\begin{itemize} 
\item Let $N_h \subset M_h$ be a component of the inverse image of $N$. 
The inverse image of $P$ in $N_h$ meets the interior of any $3$-crescent in the Kuiper completion $\che N_h$ of $N_h$. 
The fundamental group of $N$ is virtually infinite-cyclic.  
\item Let $R$ be a $3$-crescent in $\clo(N_h)$ in $\che M_h$. 
Then $R$ is a bihedral $3$-crescent and 
$\dev_h|\Lambda_1(R)$ is a homeomorphism to $H -K$ for a properly 
convex compact domain $K$ in a $3$-hemisphere $H$
with $K \cap \partial H \ne \emp$. 
\end{itemize}
\end{itemize} 
\end{theorem} 

A toral $\pi$-submanifold is {\em maximal} if no toral $\pi$-submanifold of type II 
contains it properly. 

To prove Theorems \ref{thm:caffI} and \ref{thm:caffII}, we just need to study the case when $N$ is 
a \hyperlink{term-camI}{concave affine $3$-manifold} with compressible boundary.
Let $N_{h}$ denote a component of the inverse image of $N$ in $M_{h}$ as in the premise. 

Suppose that we obtain a bihedral $3$-crescent $R$ in $\che N_h$
so that a deck transformation $g$ acts on $R^o \cup I_R^o - \{x\} \subset N_h$ properly. 
We call such a bihedral $3$-crescent a \hypertarget{term-tbc}{{\em toral bihedral $3$-crescent}}
and $g$ the \hypertarget{term-ass}{{\em associated}} deck transformation group. 



This proof is fairly long. To outline, we give the following: 
\begin{itemize} 
\item[(I)] Concave affine $3$-manifolds of type I. 
\item[(II)] Concave affine $3$-manifolds of type II. 
\begin{itemize} 
\item[(A)] There exist three mutually \hyperlink{term-ovl}{overlapping} bihedral $3$-crescents in $\Lambda(R)$ for 
a bihedral $3$-crescent $R$ in $\che M_h$. 
\begin{itemize} 
\item[(i)] There is a pair of \hyperlink{term-opp}{opposite} bihedral $3$-crescents in $\Lambda(R)$. 
By Lemma \ref{lem:casei}, $M$ is covered by an affine Hopf $3$-manifold finitely. 
\item[(ii)] Otherwise, $\dev_h| \Lambda_{1}(R)$ is a homeomorphism to $H- K$ for a properly convex domain $K$ and a $3$-hemisphere $H$
containing $K$, and  
$\Lambda(R)$ contains a toral bihedral $3$-crescent. 
Lemma \ref{lem:toralpi} gives us a toral $\pi$-submanifold.
\end{itemize}
\item[(B)] Otherwise, 
all bihedral $3$-crescents $R$ have $\dev_{h}(I_R)$ 
containing a fixed pair of points $q, q_-$. Then $\Lambda(R)$ is a union of segments from $q$ to $q_-$.
\begin{itemize} 
\item[(i)] A closed curve in a component $A_1$ of $\Bd \Lambda(R) \cap M_h$ bounds a disk in the union $A_{1, +}$ of lines from $q$ to $q_-$ passing $A_1$. 
Here, the situation is similar to (A)(i), and we use Lemma \ref{lem:toralpi}.
\item[(ii)] Otherwise, $A_{1, +}$ is an annulus. We show that this case does not happen. 
\end{itemize} 
\end{itemize}
\end{itemize}

\subsubsection{Case \ep{I} {\rm :}} \label{subsub:I} 
Let $N$ be a concave affine $3$-manifold of type I in $M$. 
Then $F \cap M_h$ covers $N$ for a hemispherical $3$-crescent $F$. 
Let $\Gamma_N$ denote the subgroup of $\Gamma_h$ acting on $F \cap M_h$ 
as the deck transformation group of the covering map to $N$. 

Let $\tilde A_1$ denote a compressible component of $I_F \cap M_h$. 
By Lemma \ref{lem:simplyc}, $\tilde A_1$ is not simply connected. 
By Corollary \ref{cor:cyclic}, 
\[\tilde A = I_F^o \cap M_h = I_F^o -\{x\}\] holds.  
Thus, 
\[ N_{h}=F \cap M_h = F^o \cup I_F^o -\{x\}.\] 



\subsubsection{Case \ep{II} {\rm :}} \label{subsub:II} 
Now suppose that $N$ is a concave affine $3$-manifold of type II in $M$.  
We assume that there is no hemispherical $3$-crescent in $\che M_{h}$. 
Then $N$ is covered by $\Lambda(R) \cap M_h$ for a bihedral $3$-crescent $R$ in $\che M_{h}$. 
Let $\Gamma_N$ denote the subgroup of $\Gamma_h$ acting on $\Lambda(R) \cap M_h$ 
as the deck transformation group of the covering map to $N$. 
Recall that 
\[ \dev_{h}(\Lambda(R)) \subset H, \dev_h(\delta_{\infty} \Lambda(R)) \subset \partial H\]
for a $3$-hemisphere $H \subset \SI^{3}$. (See Corollary 5.8 of \cite{psconv}.)

\begin{description} 
\item[(II)(A)] Suppose that there exist three mutually overlapping bihedral $3$-crescents $R_1, R_2,$ and $R_3$ with 
$\{I_{R_i}|i=1,2,3\}$ in general position.
\item[(II)(B)] Suppose that there exist no such triple of bihedral $3$-crescents.
\end{description}
We will defer (II)(B) to Section \ref{subsub:B}.  

Now assume (II)(A). 
By modifying the proofs of Lemma 11.1 and Proposition 11.1 of \cite{rdsv} 
for bihedral $3$-crescents, which are not necessarily radial as in the paper, 
we obtain that 
\begin{align} 
\dev_h: & \Lambda_1(R)  \ra H  - K, \label{eqn:Lhom} \\
\dev_h| \Lambda_1(R) \cap M_h:  & \Lambda_{1}(R) \cap M_h \ra H^{o} -K \nonumber
\end{align} 
are homeomorphisms for a $3$-hemisphere $H$ and a nonempty 
compact properly convex set $K$.
(For Lemma 11.1 and Proposition 11.1 of \cite{rdsv}, we do not need $I_{R}$ for each bihedral $3$-crescent to contain the origin.
There is a mistake in the third line of the proof of Lemma 11.1 of \cite{rdsv}. We need to change 
$P_{1} \cap L_{1}$ and $P_{1}\cap L_{2}$ to $P_{1}\cap L_{2}$ and $P_{1}\cap L_{3}$ respectively. )
The general position property of $I_{R_i}$, $i=1, 2, 3$, implies that $K$ is properly convex.
Also, \eqref{eqn:Lhom} implies that $h_h|\Gamma_N$ is injective. 

Here, $\dev_h(\alpha_{R'}) \subset \partial H$ for $R' \sim R$.  
Now, $\Bd \Lambda_1(R) \cap M_h$ is mapped into $\Bd K$. 
Let $K'$ denote the inverse image in $\Bd \Lambda_1(R)$ of $K$. 
$h_h(\Gamma_N)$ is an affine transformation group of $H^o$ since it acts on 
an affine space $H^o$ as a projective automorphism group. 

\begin{hyp}
We can have two possibilities: 
\begin{description}
\item[(II)(A)(i)] Suppose that there exist two opposite bihedral $3$-crescents $S_1, S_2 \sim R$. 
\item[(II)(A)(ii)] There are no such bihedral $3$-crescents. 
\end{description} 
\end{hyp}



\subsubsection{Case\ep{II}\ep{A}\ep{i} {\rm :}}  \label{subsub:IIAi} 
At least one component $A_1$ of $I_{S_1}^o \cap M_h$ contains $I_{S_1}^o - K'$ for a properly convex compact set $K'$
by \eqref{eqn:Lhom}. 
Here, we will show that $M$ is an affine Hopf $3$-manifold. 
The following finishes (A)(i). 
\begin{lemma} \label{lem:casei} 
Suppose that there exist two bihedral $3$-crescents $S_1, S_2 $ in $\che M_h$ so that $I_{S_1} \cap M_h$ and $I_{S_2} \cap M_h$ 
intersect and are tangent but $\dev_h(S_1)^o \cap \dev_h(S_2)^o =\emp$. 
Assume $S_1, S_2 \sim R$, and \ep{II}\ep{A}\ep{i}. 
Then there exists a unique component of $I_{S_i} \cap M_h$ equal to $I_{S_i}^o-\{x\}$ for a point $x$ of $I_{S_i}^o$, $i=1, 2$,
and $M$ is an affine Hopf $3$-manifold.
\end{lemma} 
\begin{proof} 
First, $I_{S_1} \cap M_h$ and $I_{S_2} \cap M_h$ meet at the union of  their common components by geometry
since such a component is totally geodesic and complete in $M_h$.

From above, $A_{1}$ is a quotient of a domain  containing $I_{S_1}^o - K'$ for a properly convex compact set $K'$. 
By the classification of the affine $2$-manifolds (see \cite{BenTor}), the only possibility is 
\[ A_{1} = \begin{cases} I_{S_1}^o \hbox{ or } \\I_{S_1}^o -\{x\}, x \in I_{S_1}^o. \end{cases}\] 

In the first case, 
we obtain that \[\dev_h(\Lambda_1(R) \cap M_h)=H^o \hbox{ and }\]
\[\partial H = \clo(\alpha_{S_{1}}) \cup \clo(\alpha_{S_{2}}) \subset \che M_{h,\infty}.\] 
Hence, $M_h$ is projectively diffeomorphic 
to the complete affine space. Thus, $M_h$ is diffeomorphic to $\bR^3$, a contradiction to the assumption. 

Now suppose that $A_1 = I_{S_1}^o - \{x\}$. 
Since \[\clo(\alpha_{S_1}) \cup \clo(\alpha_{S_2})\subset M_{h, \infty},\] 
$S_1^o \cup A_1 \cup S_2^o$ is homeomorphic to $\SI^2 \times \bR$ 
and $M_h = S_1^o \cup A_1 \cup S_2^o$.
Since $S_{1}$ and $S_{2}$ are mapped into the closures of two different components of $\SI^{3} - \SI^{2}$ respectively,
\[\dev_{h}|S_1^o \cup A_1 \cup S_2^o\] is an embedding onto 
its image by geometry. 
Since $\dev_{h}(x)$ is an isolated boundary point, 
Corollary \ref{cor:cyclic} implies the result that $M$ is an affine Hopf manifold. 
\end{proof}


\subsubsection{Case \ep{II}\ep{A}\ep{ii} {\rm :}}  \label{subsub:IIAii} 
In this case, $K^o$ is a nonempty properly convex open domain, and $\Bd \Lambda(R) \cap M_h$ is mapped into $\Bd K$: 
Otherwise, $\dim K \leq n-1$ and $K$ is a subset of a hyperspace $V$.
Then the two components of $H^o- V$
lift to cell embeddings in $\Lambda(R)^o$ by \eqref{eqn:Lhom}.  The closures of two cells in $\Lambda(R)$ are 
bihedral $3$-crescents again by \eqref{eqn:Lhom}. 
The two crescents are opposite. Thus, we are in case (i), a contradiction. 


By Lemma \ref{lem:BdL1BdL} and \eqref{eqn:Lhom}, we have
$\Bd \Lambda_1(R) \cap M_h=\Bd \Lambda(R) \cap M_h$. 
For following Lemma \ref{lem:BdL1BdL}, we don't need to assume (II)(A) 
and $K$ needs not be properly convex. 
\begin{lemma} \label{lem:BdL1BdL} 
Suppose that 
$\dev_{h}(\Lambda_{1}(R))$ is a homeomorphism to $H - K$ for  a convex domain $K^{o}\ne \emp$. 
Then $\Bd \Lambda_1(R) \cap M_h=\Bd \Lambda(R) \cap M_h$. 
\end{lemma} 
\begin{proof} 
Since for each crescent $S$, $S^{o}$ is dense in $S$, we obtain
\[\Bd \Lambda(R) \cap M_{h}\subset \Bd \Lambda_{1}(R)\cap M_{h}.\] 
Given a point $x \in \Bd \Lambda_{1}(R)$, choose a convex open neighborhood $B(x) \subset M_{h}$
with $\dev_{h}|B(x)$ is an embedding. 
$B(x) \cap S^{o}$ for a crescent $S$, $S \sim R$, is 
a closure of a component of $B(x) - I_{S} \cap B(x)$ for a totally disk $I_{S}\cap B(x)$ with boundary in $B(x)$. 
The set $B(x) - \Lambda_{1}(R)$ is a convex set $K''$ in $B(x)$.
Since $\dev_{h}(\Lambda_{1}(R))$ is a homeomorphism to $H - K$ by premise, 
\begin{itemize}
\item $\dev_{h}(x) \in K$, and 
\item $\dev_{h}| B(x) \cap \Lambda_{1}(R)$ is an embedding to $\dev_{h}(B(x)) - K$, and hence, 
\item $\dev_{h}| B(x) - \Lambda_{1}(R) (= K'')$ is an embedding to  $\dev_{h}(B(x)) \cap K$. 
\end{itemize}

Suppose that $K''$ has an empty interior.  Then  
$\dev_{h}(K'')$ has an empty interior. This implies by geometry of the supporting 
hyperplanes that $K^o$ empty which we showed to be absurd above
in this subsection. 

Now $K$ has a nonempty interior. 
The interior of $K$ is disjoint from $\dev_{h}(T)$ for any crescent $T$, $T \sim R$ since otherwise 
\[\dev_{h}(T^{o}) \cap K^{o} \ne \emp \hbox{ while } K^{o} \cap \dev_{h}(\Lambda_{1}(R))= \emp.\]  
Thus, $K^{\prime \prime o} \cap \Lambda(R) = \emp$, and $x \in \Bd \Lambda(R)$. 
\end{proof}

\begin{lemma}\label{lem:Kncompact}
Assume as in Theorem \ref{thm:caffII} with a concave affine $3$-manifold $N$ with compressible boundary. 
Suppose that $N$ is in case {\rm (II)(A)(ii)}. Let $K$ be the complement of $\dev_{h}(\Lambda_{1}(R))$. 
Then $K$ is an unbounded subset of an affine space $H^o$. 
Moreover, $K \cap \partial H$ is a nonempty compact convex set, 
and $\Bd K \cap H^{o}$ is homeomorphic to a disk. 
\end{lemma} 
\begin{proof} 
Suppose that $K$ is a bounded subset of $H^{o}$. Then $\dev_h| \Lambda_{1}(R)$ is 
a homeomorphism to $H - K$ by \eqref{eqn:Lhom}. 
By Theorem \ref{thm:bdsph}, components of $\Bd \Lambda(R) \cap M_{h} = \Bd \Lambda_{1} \cap M_{h}$ by Lemma \ref{lem:BdL1BdL} 
are not a sphere. 
$\Bd \Lambda_{1}\cap M_{h}$ is mapped into a surface in $\partial K$.  Since the map cannot be onto the sphere $\partial K$, 
there exists a noncompact component $A_1$ of $\Bd \Lambda(R) \cap M_h = \Bd \Lambda_{1} \cap M_{h}$ covering 
a closed surface $B_1$. 
By Lemma \ref{lem:nobd}, this is a contradiction as $h_h| \Gamma_N$ is injective by \eqref{eqn:Lhom}. 
Hence, $K$ is unbounded in $H^{o}$. 
\qed
\end{proof}





\begin{figure}[h]

\begin{center}
\includegraphics[height=6cm]{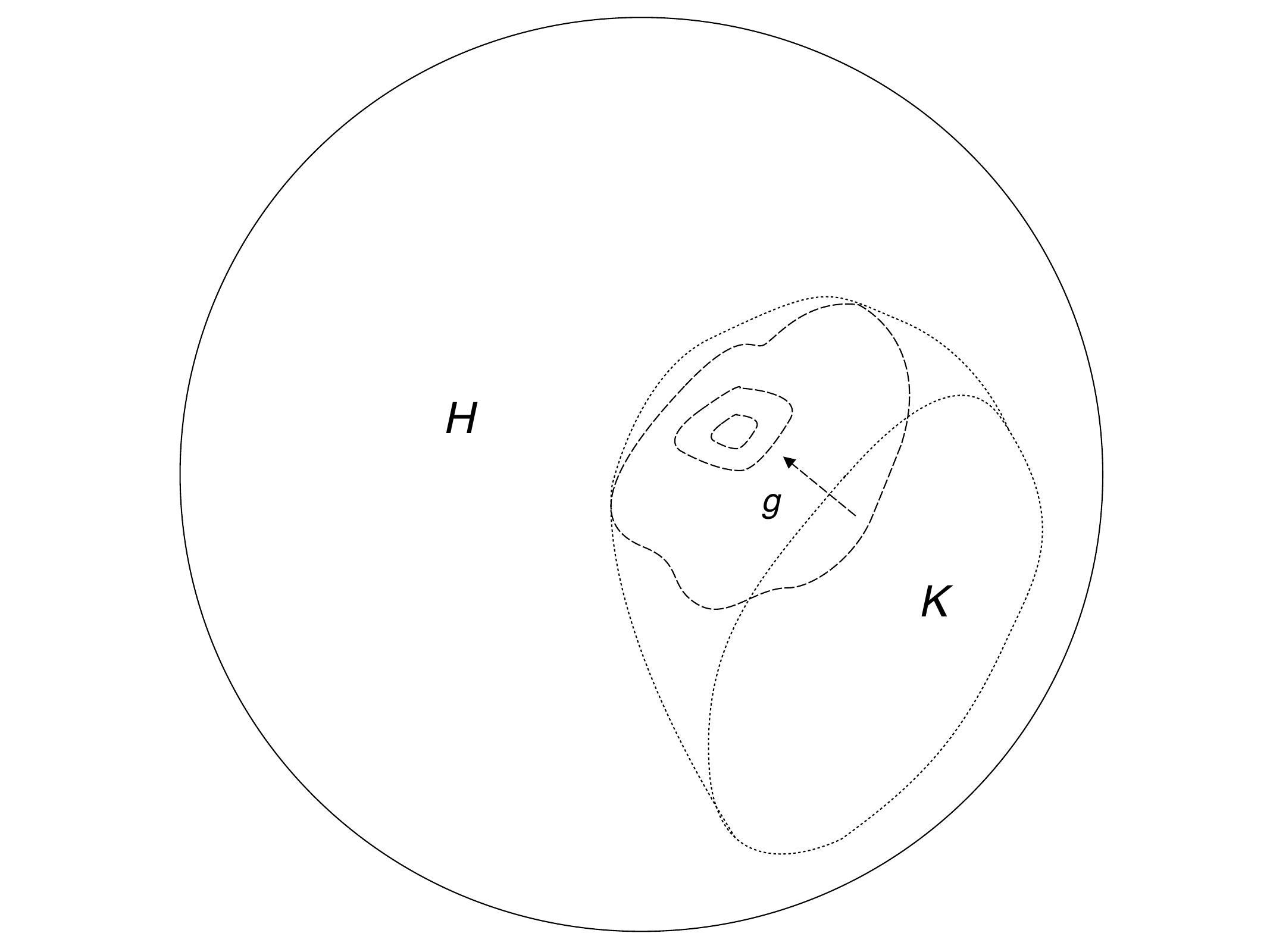}
\end{center} 
\caption{The diagram for $H, K$ and $c_{1}\subset \Bd K$ for the case (II)(A)(ii).} 
\label{fig:concave}
\end{figure}

\subsubsection{Case \ep{II}\ep{A}\ep{ii} {\rm :} obtaining a toral bihedral $3$-crescent} \label{subsub:IIAii-2}
The next step is to show that $\dev_h| \Lambda(R): \Lambda(R) \ra H$ is injective: 

If each component of $\Bd \Lambda(R) \cap M_h $ is simply connected, then 
each component is a disk by Theorem \ref{thm:bdsph}.
By Lemma \ref{lem:simplyc}, these components are $\pi$-injective. 
Hence, a component $A_1$ of $\Bd \Lambda(R) \cap M_h$ is not simply connected 
and becomes compressible in $N$. 

Now, $\clo(K) \cap \partial H = \Bd K \cap \partial H$ is a convex compact disk in $\Bd K$, 
and $\Bd K$ is homeomorphic to a sphere. 
Since \[\Bd K \cap H^{o } = \Bd K - \Bd K\cap \partial H\] holds, 
$\Bd K \cap H^{o}$ is homeomorphic to a disk as surface topology tell us. 
Since $A_1$ is not simply connected, 
there exists a simple closed curve $c_1 \subset A_1$ so that $\dev_{h}(c_1)$ bounding a disk $D_1$ in $\Bd K \cap H^{o}$. 
By Lemmas \ref{lem:diskfixedpt} and \ref{lem:toralpi}, $N$ contains a toral $\pi$-submanifold. 
This finishes the case (II)(A)(ii). 

\begin{lemma}\label{lem:diskfixedpt}  
Assume	 \eqref{eqn:Lhom}, and 
suppose that  $\Bd \Lambda(R) \cap M_{h}=\Bd \Lambda_{1}(R)\cap M_{h} $, and  it is mapped
into a properly embedded disk $S$ in $H^{o}$ bounding a convex domain $K$ in 
the complement of $\dev_{h}( \Lambda(R) \cap M_{h})$. 
Let $A_{1}$ be a non-simply connected domain in 
$\Bd \Lambda(R) \cap M_{h}$ containing a curve $c_{1}$ so that $\dev_{h}(c_1)$ bounding an open disk $D_1$ in $S$. 
Assume that $\Gamma_{N}$ acts on $S$. 
Then the following hold {\rm :} 
\begin{itemize} 
\item $\dev_{h}| \Lambda(R): \Lambda(R) \ra H$ is injective. 
\item $\Bd \Lambda(R) \cap M_{h}$ has a unique component.
\item $\Lambda(R)$ contains a toral bihedral $3$-crescent $R_{P}$. 
\item The fundamental group of $N$ is virtually infinite cyclic. 
\end{itemize} 
\end{lemma}
\begin{proof} 
Let $\Gamma_1$ be the subgroup of $\Gamma_N$ acting on $A_1$ cocompactly. 
There exists an element $g \in \Gamma_1$ such that 
$h_h(g)(\dev_h(c_1)) \subset D_1 \cap \dev_h(A_1)$ by Lemma \ref{lem:disj}. 
Also, $h_{h}(g)(D_{1}) \subset D_{1}$ since the outside component of $S - \dev_{h}(c_{1})$ is not homeomorphic to a disk. 
By \eqref{eqn:Lhom}, $\Lambda_1(R)$ is homeomorphic to $H^o - K$, a cell. 
We find an open disk $D'_1$ in $\Lambda_1(R)$ that is compactified by boundary $c_1$. 
Then $D'_1$ bounds a $3$-dimensional domain $B_1$ in $\Lambda_1(R) \cap N_h$ where $B_{1}^{o}$ is a cell. 
Since the group action is proper, and the closure of $D'_{1}$ in $N_{h}$ is compact, 
we may also assume that $g(B_1) \subset B_1^o$ disjoint from $D'_1$ by Lemma \ref{lem:disj}. 
Thus, we can find a fixed point $x$ in $\clo(K)$ for $g \in \Gamma_1$ by the Brouwer fixed-point theorem.
We can verify the premises of Proposition \ref{prop:diskf} for $\dev_h(D_{1}')$ 
by using a supporting hyperplane at $x$ since $K$ is convex and the boundary of $\dev_h(D_1')$ is in $K$. 
Proposition \ref{prop:diskf} implies that 
$x$ is the fixed point of the largest norm eigenvalue of 
$h_h(g)$ and the global attracting fixed point of $h_{h}(g)| H^{o}$. 




Now, we prove the injectivity of $\dev_{h}| \Lambda(R)$: 
Let $x_{j}, j=1, 2$ be points of $\Lambda(R)$.
Let $R_1, R_2 \sim R$ be two bihedral $3$-crescents where $x_{j} \in R_{j}, j=1, 2$. 
We may assume that $\dev_{h}(R_j)$ meets $\dev_{h}(\Bd \Lambda_1(R)) - \partial H$ by 
taking the maximal bihedral $3$-crescents. 
Then $h_{h}(g)^i\dev_{h}(R_j)$ meets a neighborhood of $x$ for sufficiently large $i$
by \eqref{eqn:Lhom}. 
Since $\dev_{h}(g^{i}(I_{R_{1}}))$ and $\dev_{h}(g^i(I_{R_{2}}))$ are very close containing nearby points for sufficiently large $i$
and supporting a properly convex domain $K^o$,  we obtain that  $\dev_h(g^i(R_1))$ and $\dev_h(g^i(R_2))$ meet in the interior. 
By \eqref{eqn:Lhom}, we obtain \eqref{eqn:Lhom},
\[g^i(R_1)^o \cap g^i(R_2) \ne \emp\]{ for sufficiently large $i$ and hence }
\[ R_1^o \cap R_2^o \ne \emp.\]
By Theorem 5.4 and Proposition 3.9 of \cite{psconv}, 
$\dev_h|R_1\cup R_2$ is injective.  
Therefore, $\dev_h|\Lambda(R)$ is injective.  
This proves the first item.

Since $\dev_h|\Lambda(R) \cap M_h$ is injective,  the restriction of an immersion 
\[\dev_h|\clo(K) \cap \Bd \Lambda(R) \cap M_h\] is a homeomorphism to its image $Y$ in 
$\Bd K$. The set $Y$ is an open surface. Then $Y/h_h(\Gamma_N)$ is a union of closed surfaces. 
Let $Y_1$ be the image of $A_1$. $Y_1/h_h(\Gamma_1)$ is a connected closed surface homeomorphic to $A_1/\Gamma_1$. 



Since $\dev_h(x)$ is a unique attracting fixed point of $h_h(g)$ in $H^{o}$, 
$h_h(g)^i(c_1)$ goes into an arbitrary neighborhood of $\dev_h(x)$ in $\Bd K$ for sufficiently large $i$. 
$h_h(g)^i(c_1)$ goes into an arbitrary tubular neighborhood of $\Bd K \cap \partial H$ in $\Bd K$ for sufficiently small
negative number $i$. Using $i$ and $-i$ for a large integer $i$, $h_h(g)^i(c_1)$ and $h_h(g)^{-i}(c_1)$ bound a compact annulus in $\Bd K \cap H^o$. 
If there is a component $\tilde Y_j$ of $Y \subset \Bd K$ other than $\dev_h(A_1)$, then it lies in one of the annuli, a bounded subset of $H^o$,  
and $\tilde Y_j$ covers a compact surface $Y_j$ for some $j$. 
By Lemma \ref{lem:nobd}, this is a contradiction. 
Thus, $\Bd \Lambda(R) \cap M_h$ has a unique component. 

Let $D$ denote the disk $\Bd K \cap H^o$.
$Q_K:= D - \dev_h(\tilde A_1)$ is either 
$\{\dev_h(x)\}$ or a closed set with infinitely many 
components since $h_h(g)^i(c_1)$ are disjoint from $Q_K$ and 
$\langle h_h(g)\rangle$ acts on $Q_K$. 

We obtain $Q_K =\{\dev_h(x)\}$ by Lemma \ref{lem:fullset}
and
\begin{equation}\label{eqn:A1bdK}
(\Bd K -\{\dev_h(x)\}) \cap H^o = \dev_h(A_1).
\end{equation}

Since $\Gamma_N$ acts faithfully, properly discontinuously, and freely on an annulus $A_1$,
$\Gamma_N$ is virtually infinite-cyclic. 
The existence of $g$ shows that the $h_h(\Gamma_N)$ fixes the unique point $\dev_h(x)$ 
corresponding to one of the ends. 
This proves the fourth item. 


Let $K_x \subset \SI^{2}_x$ denote the subspace of directions of the segments with endpoints in $\dev_{h}(x)$ and $K^o$. 
Obviously, $K_x$ is a convex open domain in an open half-space of $\SI^{2}_x$. 
Our $h_{h}(g)$ acts on $K_x$, and 
$\SI^{2}_x$ has a $h_{h}(g)$-invariant great circle $\SI^1$ outside $K_x$
as we can  deduce by the existence of $K_{x}$.

We take a union of maximal segments in $\dev_{h}(\Lambda(R))$ from $\dev_{h}(x)$ 
in directions in $\SI^{1}$. Their union is 
a $2$-hemisphere $P$ with boundary in $\partial H$, and $\dev_{h}(x) \in P$.


We find an open bihedron $B \subset H-K$ whose boundary contains 
an open $2$-hemisphere in $\partial H$ and $P$. 
By taking an inverse and the closure, we obtain a bihedral $3$-crescent $R_P \subset \Lambda(R)$ with
$x \in I_{R_P}$. 
By the first item, $g$ acts on $R_{P}$, $I_{P}$, and $x$.

The last step is to show $R_P$ has the desired property. 
By our choice of $K_{x}$ and $P$, we obtain 
$\dev_{h}(I_{R_{P}})^{o} - \dev_{h}(x) \subset H - K^{o}$.
By \eqref{eqn:A1bdK} and the first item, we obtain
\begin{equation} \label{eqn:A1}
I_{R_{P}}^{o} - \{x\} \subset A_{1} \cup \Lambda_1(R).
\end{equation} 
Hence, 
$I_{R_P}^o -\{x\} \subset N_h$ for our bihedral $3$-crescent $R_P$ above.
There is an element $g \in \Gamma_{N}$ acting on $R_{P}^{o}\cup I_{R_{P}}^{o} -\{ x\}$.  

\end{proof}

\begin{lemma} \label{lem:fullset} 
Let $S_{0}$ be a properly embedded disk or cylinder in $\bR^{3}$. 
Let $\tilde A_{0}\subset S_{0}$ be a connected open set covering a closed surface $A_{0}$
with the deck transformation group $G_{1}$ also acting on $S_{0}$ for $G \subset \Aff(\bR^{3})$.
Suppose that there exists a collection of simply closed curves $c_{i} \in A_{0}$, $i \in \bZ$, 
so that for any end neighborhood of $S_{0}$ there is a component of $S_{0} - c_{i}$ in it. 
Then $S_{0} - \tilde A_{0}$ cannot have infinitely many components. 
\end{lemma} 
\begin{proof} 
Suppose that $\tilde A_{0}$ is an open planar surface with infinitely many ends. 
Giving an arbitrary complex structure on $A_{0}$, the cover 
$\tilde A_{0}$ admits a Koebe uniformization as 
$\bC P^1 - \Lambda$ for a Cantor set $\Lambda$. 
That is $A_{0} = \tilde A_{0}/\Gamma_1$ is homeomorphic to a closed Schottky Riemann surface 
$(\bC P^1 - \Lambda)/G_1$ where $G_1$ is in $\PSL(2, \bC)$ isomorphic to $\Gamma_1$. 
(See p. 77 of Marden \cite{Marden} for the proof.)
The set of the pairs of fixed points of elements of $G_1$ 
are dense in $\Lambda \times \Lambda - \Delta(\Lambda)$ 
for the diagonal $\Delta(\Lambda)$ of $\Lambda$. 
(See Theorem 2.14 of Apanasov \cite{Apanasov}.) 
We can find a closed curve $c$ 
in the surface $\tilde A_{0}/\Gamma_1$ so that $c$ lifts to 
a curve $\tilde c$ ending at two points in arbitrary two small open neighborhood of 
two points in $\Lambda$. Let $k_1$ and $k_2$ be two points of $\Lambda$. 
On $\tilde A_{0}$, simple closed curves bound end neighborhoods. 
We may assume that $k_1$ corresponds to an end of $\tilde A_{0}$ whose
end neighborhood is bounded by a simple closed curve $d_1$, 
and $k_2$ corresponds to an end of $\tilde A_{0}$ whose end neighborhood
is bounded by a simple closed curve $d_2$. 
We assume that $\dev_h(d_i)$ bounds an open disk $D_i$ for $i=1, 2$, whose closure is 
compact in $D$ since we can choose $k_1$ and $k_2$ arbitrarily. 
Assume that $D_1 \cap D_2 = \emptyset$.

Choose an orientation-preserving 
element $g_c \in \Gamma_1$ acting on $\tilde c$.
Then \[h_h(g_c^n)(\dev_h(d_1)) \subset D_2 \hbox{ for some } n.\] 
By orientation considerations of how $\tilde c$ meets 
$d_1$ and $g_c^n(d_1)$, we obtain
\[h_h(g_c^n)(D - D_1) \subset D_2.\] 
Since $g_c$ acts on $\partial H$, $D - D_1$ has limit points in 
$\partial H$, and $D_2$ has no limit points in $\partial H$, 
this is a contradiction. 
\end{proof}

\subsubsection{Case \ep{B} {\rm :}} \label{subsub:B} 
Now suppose that $\Lambda(R)$ contains no triple of 
mutually \hyperlink{term-ovl}{overlapping} bihedral $3$-crescents $S_{i}$, $i=1,2,3$, with $\dev_{h}(I_{S_{i}})$ in general position.
By induction on overlapping pairs of bihedral $3$-crescents, we obtain that 
$\dev_h(I_S)$ for a bihedral $3$-crescent $S, S \sim R,$ share a common point $q \in \partial H$ and hence its antipode $q_- \in \partial H$. 
Then $\Lambda(R)$ is a union of segments whose developing images end at $q, q_-$. 
The interior of such segments in $\Lambda(R)$ is called a {\em complete $q$-line}. 
Also, {\em $q$-lines} are subarcs of complete $q$-lines. 
$\Bd \Lambda(R) \cap M_h$ is foliated by subsets of $q$-lines. 



Let $\SI^2_q$ denote the sphere of directions of complete affine lines from $q$,
and let $\SI^2_q$ have a standard Riemannian metric of curvature $1$. 
The space of $q$-lines in $\Lambda_1(R)\cap M_h$ whose developing image go
from $q$ to its antipode $q_-$  is an open surface $S_R$ with an affine structure. 
The developing map $\dev_h$ induces an immersion $\dev_{h, q}: S_R \ra \SI^2_q$.
The surface $S_R$ develops into a $2$-hemisphere \hypertarget{not-hq}{$H_q \subset \SI^{2}_{q}$ }
whose interior $H_q^o$ is identifiable with an affine $2$-space. 
Denote by $\Pi_{q}: H^{o} \ra H^{o}_{q}$ the projection. 

The \hyperlink{term-kcp}{Kuiper completion} $\che S_R$ of $S_{R}$ has an ideal subset $c'$ 
that is the image of $\Bd \Lambda(R) \cap M_{h}$ 
and a geodesic boundary subset corresponding to $\delta_\infty \Lambda(R)$ and is mapped to $\partial H_q$.
We denote the extension by the same symbol $\dev_{h, q}: \che S_{R} \ra \SI^{2}_{q}$.



Denote $N_h := \Lambda(R) \cap M_h$ covering a concave affine manifold in $M_h$. 
If each component of $\partial N_h = \Bd \Lambda(R) \cap M_h$ is simply connected, 
then it is incompressible by Theorem \ref{thm:bdsph} as in the beginning of Section \ref{subsub:IIAii-2}.
Thus, there is a component $A_{1}$ 
of $\Bd \Lambda(R) \cap M_h$ containing a simple closed curve $c$ 
that is not null-homotopic in $A_{1}$. 
We will use the same notation $\dev_h$ for the extension of $\dev_h|N_{h}$ to $\che N_h$. 
Let $\mathcal{L}_{q}$ denote the set of complete $q$-lines $l$ such that 
\[l \subset R''' \hbox{ for } R''' \sim R, R''' \subset \che N_{h}, \hbox{ and } l \cap A_{1} \ne \emp.\] 
We define 
\[ A_{1+}:= \bigcup_{l \in \mathcal{L}_{q}} l.\]

We claim that $A_{1+}$ is homeomorphic to the injective image of a topologically open surface: 
Recalling the surface $S_R$ above, we obtain a fibration $\Pi_{R}: \Lambda_1(R) \cap M_h \ra S_R$
extending to $\che N_{h} \ra \che S_{R}$, to be denoted by $\Pi_{R}$ again.  
$\Pi_{R}$ maps $A_{1+}$ to a set $a_{1+}$ in the \hyperlink{term-ids}{ideal boundary} of $\che S_R$ of $S_R$. 
Since $q$-complete lines pass the open surface $A_{1}$ 
foliated by $q$-arcs and 
$\dev_{h,q}|a_{1+}$ maps locally injectively to an embedded arc in \hyperlink{not-hq}{$H^{o}_{q}$}. 
Thus, $a_{1+}$ is a locally injective open arc
since $A_1$ is a surface. 

Suppose that two leaves $l_1$ and $l_2$ of $A_{1+}$ 
go to the same point of an open arc $\alpha$ in $a_{1+}$ where $\dev_{h, q}| \alpha$ is an embedding.  
Since $l_{1}$ and $l_{2}$ are fibers, there is a point $\Pi_R(l)$ in $S_{R}$ of $\bdd$-distance $< \eps$ from the images 
$\Pi_R(l_{1}), \Pi_R(l_{2})$ of these lines in $\che S_{R}$. 
Inside $\Lambda_1(R)$, there exist paths of $\bdd$-length $< \eps$ from $l_1$ and $l_2$ to any point of 
a common line $l$ in $\Lambda_1(R)$ corresponding to $\Pi_R(l)$ by spherical geometry. 
Taking $\eps \ra 0$ and $l$ closer to $l_i$, we obtain $l_1=l_2$. 
Hence, we showed that $A_{1+}$ fibers over $a_{1+}$ locally. 

This implies that $A_{1+}$ is the image of an open surface.
We give a new topology on $A_{1+}$ by giving a basis of $A_{1+}$ as the set of components of 
the inverse images of open sets in $\che M_h$. Then $A_{1+}$ is homeomorphic to  a surface with this topology.



As above, $A_1$ contains a simple closed curve $c$ not homotopic to a point in $A_1$. 
\[c''':= \Pi_{R}(c) \subset a_{1+}\] 
is either a compact arc, i.e., homeomorphic to an interval or 
a circle. 

We divide into two cases: 
\begin{description} 
	\item[(i)] $c'''$ is homeomorphic to an interval. 
	\item[(ii)] $c'''$ is homeomorphic to a circle. 
\end{description}

(i) 
Then $c$ bounds an open disk $D'$ in the fibered space $A_{1+}$. 
Let $\Gamma_1$ be the subgroup of $\Gamma_N$ acting on $A_1$.
We can use a similar argument to (II)(A)(ii):
First, there exists $g \in \Gamma_1$ so that $g(c)$ is in $D' \cap A_{1}$ by Lemma \ref{lem:disj}. 
Hence $g$ fixes a point $x$ in $D^{\prime o}$ that is a fixed point on $A_{1+}$
by the Brouwer fixed-point theorem.  
$A_{1+}$ is either homeomorphic to an annulus or a disk since $A_{1+}$ is foliated by $q$-lines.
We have $g(D') \subset D'$ since 
exactly one component of $A_{1+} - g(c)$ is homeomorphic to a disk.




Let $x_{q}$ denote  the complete $q$-line containing $x$ in $a_{1+}$.
Let $g_{q}: \che S_{R} \ra \che S_{R}$ be the induced map of $g: \Lambda(R) \ra \Lambda(R)$. 
Recall the affine space \hyperlink{not-hq}{$H^{o}_{q}$}. 
Consider $x_{q}$ as the origin. 
Since the induced linear transformation $h(g)_{q}:H^{o}_{q} \ra H^{o}_{q}$ is not trivial, 
$h(g)_{q}$ has an isolated fixed point or has a line $l$ of fixed points in
$H^o_q$. In the second case, $h(g)_q$ acts on lines parallel to $l$, or 
acts on a parallel set of lines transversal to the line $l$.
%
The action on $\dev_{h, q}(a_{1+})$ of $h(g)_q$,
its fixed point $x_{q}$ in $H^{o}_{q}$ is locally isolated, 
or there is a geodesic subarc of fixed points forming a neighborhood or a one-sided 
neighborhood of $x_q$ in the local arc $\dev_{h, q}(a_{1+})$.



We consider the first case.  
Since $g(c)$ is in the open disk in $A_{1+}$ bounded by $c$, 
a compact arc neighborhood of 
$x_q$ in $a_{1+}$ goes into itself strictly under $g_q$. 
It must be that 
$\dev_{h, q}(x_q)$ is the attracting fixed point under $h(g)_q$. 
Thus, the local arc $\dev_{h, q}(a_{1+})$ is the union 
$\bigcup_{i \geq 0} h(g)_q^i(I)$ for a small embedded open arc $I$ 
containing $\dev_{h, q}(x_q)$.
Since $I$ is embedded, 
$\dev_{h, q}(a_{1+})$ is also an embedded arc. 
By the classification of the infinite cyclic linear automorphism groups of $H^{o}_{q}$, 
we can show that $\dev_{h, q}(a_{1+})$ is a properly embedded convex arc in $H^{o}_{q}$. In the other cases, a similar argument will show 
the same fact where we replace $I$ with an $\epsilon$-$\bdd$-neighborhood of 
the straight arc neighborhoods.

Consider the commutative diagram
\begin{alignat}{3}
A_{1+}                \quad \quad   & \stackrel{\Pi_{R}}{\longrightarrow} & a_{1+}  \nonumber \\
\downarrow \dev_{h}   &        &\quad \dev_{h, q} \downarrow \nonumber  \\ 
H^{o}       \quad \quad   & \stackrel{\Pi_{q}}{\longrightarrow}    &  H^{o}_{q}.
\end{alignat}
Since the left arrows of the above commutative diagrams are fibrations, 
\[\dev_{h}| A_{1+}: A_{1+} \ra H^o\]
is a proper embedding to $H^{o}$. 
Since $\dev_{h, q}(a_{1+})$ is a properly embedded convex arc, 
$A_{1+}$ is a properly embedded surface bounding a convex domain $K$ in $H^o$. 

We claim that $\dev_{h}: \Lambda_{1}(R) \ra H$ is an embedding: 
$\dev_{h}| R^{o} \cup \alpha_{R}$ is an embedding. 
We can choose a crescent $S$ overlapping with $R$ and $\dev_{h}(I_{S})^{o}$
containing a generically chosen $y \in A_{1}$ so that the closure of the arc 
$\dev_{h, q}(\alpha_{S})$ does not contain the endpoint of $a_{1+}$.
Then for any crescent $T$ overlapping with $S$, $\dev| S^{o}\cup \alpha_{S}\cup T^{o}\cup \alpha_{T}$ is an embedding
 by Proposition \ref{prop:transversal}. 
For any crescent $T_{1}$ overlapping with $T$, 
since $\alpha_{T_{1}}$ is not antipodal to $\alpha_{S}$ by our choice, 
and $\dev_{h}| T^{o}\cup \alpha_{T} \cup T_{1}^{o} \cup \alpha_{T_{1}}$ is injective, 
$T_{1}$ overlaps with $S$ also. This implies that 
\[\dev_{h}| S^{o} \cup \alpha_{S} \cup T^{o} \cup \alpha_{T} \cup T_{1}^{o}\cup \alpha_{T_{1}}\] 
is injective. 
By induction, we obtain that 
$\dev_{h}| \Lambda_{1}(R) $ is an embedding into $H$. 

A {\em semi-affine-plane} is the closure of a component of an affine plane with a complete line removed. 
Suppose that $J:=H - \dev_{h}(\Lambda_{1}(R))$ has empty interior. 
Since $J$ is in the complement of some image of the crescents, 
$J$ is the closure of a semi-affine-plane. 
Since each point of $A_{1}$ is in a crescent $S$, $S\sim R$ with $S^{o}\subset \Lambda_{1}(R)$, 
we have $A_{1}\subset \Bd \Lambda_{1}(R)$. 
Hence $\dev_{h}(A_{1})$ is in $J$, and hence so is $\dev_{h}(A_{1+})$. 
Since $a_{1+}$ is a properly embedded open arc, $A_{1+}$ is a complete affine plane, a contradiction. 

Since $J^{o}\ne \emp$, we obtain
 $\Bd \Lambda_{1}(R) \cap M_{h}= \Bd \Lambda(R)\cap M_{h}$ by Lemma \ref{lem:BdL1BdL}. 
By Lemmas \ref{lem:diskfixedpt} and \ref{lem:toralpi}, we obtain a toral $\pi$-submanifold from the bihedral $3$-crescent $T$.

(ii) 
This case does not occur; we show that $\Lambda(R)$ is not maximal here.
Then the open surface $A_{1+}$ is homeomorphic to an annulus 
foliated by complete affine lines. 
Here, $c$ is an essential simple closed curve. 
There exists an element $g\in \Gamma_1$ sending $c$ into a component $U_{1}$ of $A_ 1 -c$ by Lemma \ref{lem:disj}. 
Replacing $U_{1}$ by $g^{i}(U_{1})$ if necessary, we may assume that 
$g(U_{1}) \subset U_{1}$. Then $g$ is of infinite order. 



We can embed $c'''$ into $\che S_{R}$.
Recall a fibration 
\begin{equation}\label{eqn:qline} 
l \ra \Lambda_{1}(R) \cap M_{h} \stackrel{\Pi_{R}}{\ra} S_{R}
\end{equation}
where fibers are $q$-lines. 
\begin{lemma} \label{lem:SRfoliate} 
$S_{R} = \Lambda(S) \cap S_{R}$ for a $2$-dimensional crescent $S$ in $\che S_{R}$. 
Also, $\che S_{R}$ is homeomorphic to a compact annulus with a boundary component $c'''$ and 
a closed curve in $\delta_{\infty}\Lambda(R)$. 
Furthermore, $\dev_{h}| \Lambda_{1}(R) \cap M_{h}$ is finite-to-one to its image. 
\end{lemma}
\begin{proof} 
The map $\Pi_{R}$ sends the interiors of bihedral $3$-crescents to the interiors of $2$-dimensional crescents. 
Since $\Pi_{R}(\Lambda_{1}(R) \cap M_{h}) = S_{R}$, we obtain the equality. 

We take for each point $z$ of $c'''$ a $2$-dimensional crescent $S_{z}$ so that $z \in I_{S_{z}}$. 
Then $\bigcup_{z\in c'''} S_{z} \cap S_{R}$ is a closed subset of $\Lambda(S) \cap S_{R}$.
By perturbation of crescents $S_{z}$, it is also open in $S_{R}$. Hence, we have 
\[\bigcup_{z\in c'''} S_{z} \cap S_{R} = \Lambda(S) \cap S_{R} \hbox{ and hence }  \bigcup_{z\in c'''} S_{z} = \Lambda(R) = \che S_{R }.  \]

Recall that $\dev_{h, q}$ sends $\che S_R$ to a $2$-hemisphere $H_q$. 
Since $c$ is a compact arc, $\dev_{h, q}|c'''$ is a map to a compact arc in $H_{q}^{o}$. 
 For each $S_{z}$, we choose a segment $s_{z}$ in $S_{z}$ connecting 
$z$ to a point of $\alpha_{S_{z}}$. The complement $\che S_{R} - \bigcup_{z\in c'''} s_{z}$
is a disjoint union of open properly convex triangles with vertices in $c'''$.
We cover each triangle by maximal segments from the vertex in $c'''$. 
We can do the blow up on $c'''$ for the vertices of the triangles so that the segments are now disjoint. 
to obtain a surface $S'$ foliated by segments mapped to segments or
segments for the triangles. Since $S$ is compact, so is $\Lambda(R)$.

Since $\che S_{R} = \Lambda(S)$, it follows that $\che S_{R}$ is a compact surface with two boundary components 
$c'''$ and another simple closed curve in $\delta_{\infty} \Lambda(S)$. 
For each point $t \in \che S_{R}$, there is a neighborhood where 
$\dev_{h, q}:\che S_{R} \ra H$ restricts to a homomorphism to an open disk with possibly an embedded arc as the boundary. 
Since $\che S_{R}$ is compact, $\dev_{h,q}$ is finite-to-one. 

Hence $\dev_{h, q}| S_{R}$ is a finite-to-one map to its image in $H_{q}^{o}$
by above. Therefore, $\dev_{h}: \Lambda_{1}(R) \cap M_{h} \ra H^{o}$ is finite-to-one to its image. 
\end{proof}


Now, $\langle h_h(g)\rangle$ acts on a nontrivial closed curve $\dev_{h, q}(c''')$ bounded in an affine space $H_q^o$ of $\SI^2_q$. 
Thus, $\langle h_h(g)\rangle$ acts as an isometry group on $\SI^2_q$ with respect to a standard metric up to a choice of coordinates on $\SI^2_q$.
Let $L(g)$ denote the linear part of $h_h(g)$ considered as an affine transformation of the affine space $H^o$. 
We obtain 
\[ h_h(g) = 
\left(
\begin{array}{cc}
 L(g)  & v(g)  \\
 0  &  1
\end{array}
\right)
\]
where $v(g)$ is a $3$-vector. 
By the classification of elements of $\SL_\pm(4, \bR)$, 
the following hold: 
\begin{itemize} 
\item $L(g)$ has the direction vector $v_{q}$ to $q$ as an eigenvector,
\item $L(g)$ induces an orthogonal linear map on
$H^{o}_{q}:=\bR^{3}/\langle v_{q} \rangle$, and 
\item we obtain $g$ by post-composing with a translation in the direction of $q$.  
\end{itemize}
Thus, $v(g)$ is in the direction of $v_{q}$. 

Suppose that $L(g)$ is parabolic with eigenvalues all equal to $1$. 
By the second property above, $L(g)$ acts as the identity on $H^{o}_{q}$, and 
$g$ is a translation on each $q$-line in $A_{1+}$. 
Since $\langle g\rangle$ acts on $A_{1+}$, and as the identity on $H^{o}_{q}$,  
$h_{h}(g)$ is of form 
\begin{equation} \label{eqn:gpara}
\left(
\begin{array}{cccc}
1  & \alpha  & \beta  & \gamma\\
0  & 1  & 0  & 0 \\
0  &  0 & 1  & 0 \\
0 & 0 & 0 & 1
\end{array}
\right) \hbox{ for } \alpha, \beta, \gamma \in \bR. 
\end{equation} 
If $\alpha, \beta$ are not all zero, then 
we can find a plane $P_{g}$ of fixed points given by $\alpha y + \beta z +\gamma = 0$
in $\bR^{3}$. The inverse image $P'_{g}$ of $P_{g}$ in $\Lambda_{1}(R) \cap M_{h}$ is not empty. 
Also, $\dev_{h}: P'_{g} \ra P_{g}$ is finite-to-one onto its image by Lemma \ref{lem:SRfoliate}. 
Then $g$ acts $P'_{g}$ with a finite order since $h_{h}(g)| P_{g}$ is the identity map. 
This contradicts our choice in the beginning of (ii). 
Since $g$ acts properly on $\Lambda_{1}(R)\cap M_{h}$, 
this is absurd. Thus, 
\begin{equation}\label{eqn:ab} 
\alpha =0, \beta=0.
\end{equation} 
Finally, $\gamma \ne 0$ since $g \ne \Idd$. Hence, 
$g$ is a translation in the direction of $v_g$. 


Otherwise, $L(g)$ acts on a one-dimensional subspace parallel to 
$v_q$ and a complementary subspace, and $h_{h}(g)$ has the form 
\[
\left(
\begin{array}{ccc}
\lambda  &  0  & \gamma\\
0  & \mu O_{g}  & 0 \\
0 & 0 & \mu
\end{array}
\right), \lambda \mu^{3} = 1, \lambda, \mu > 0
\]
for an orthogonal $2\times 2$-matrix $O_{g}$.

Now, $A_{1+}\cap M_{h}$ has a component $A_{1}$ containing $c$. 
Assume that $\mu \ne \lambda$. 
Then $g^i(c)$ geometrically converges to a compact closed curve in the interior of $A_{1+}$
as $i\ra \infty$ or $i\ra -\infty$. 
Then the limit of $\dev_h(g^i(c))$ must be on a totally geodesic subspace $P$ by the classification of elements of $\SL_\pm(4, \bR)$
passing $\dev(A_{1+})$. 
Recall \eqref{eqn:qline}. 
By Lemma \ref{lem:SRfoliate}, the annulus 
$\che S_{R}$ has $c'''$ as a boundary component. 
$S_{R}$ contains an annulus $A_{1, R}$ with boundary $c'''$. 
The inverse image $P'$ of $P$ under $\dev_{h}$ 
contains an annulus $A'_{1, R}$ embedding to $A_{1, R}$ under $\Pi_{R}$.
Then $\dev_{h}| A'_{1, R}$ is a finite-to-one map. 
We may assume that $g$ acts on $A_{1, R}$, and hence $g$ acts on $A'_{1, R}$. 
Since $g^i$ is represented as a sequence of 
uniformly bounded matrices on $A'_{1, R}$ for every $i \in \bZ$,  
and $g$ is of infinite order, this is a contradiction to the properness of the action of $\langle g\rangle$. 
Therefore, we obtain
\begin{equation}\label{eqn:ml}
\mu = \lambda = 1.
\end{equation}

Since $h_{h}(g)(q)=q$, and $h_{h}$ acts on $H^{o}$, it follows that
$h_h(g)$ restricts to an affine transformation in $H^o$ acting on the set of a parallel collection of lines. 
$h_h(g)$ acts as a translation composed with a rotation on $H^o$ with respect to a Euclidean metric
since the $3\times 3$-matrix of $L(g)$ decomposes into an orthogonal $2\times 2$-submatrix and the third diagonal element equal to $1$. 

Thus, in all cases as indicated by equations \eqref{eqn:ab} or \eqref{eqn:ml},  $g$ is of form 
\[
\left(
\begin{array}{ccc}
1 &  0  & \gamma\\
0  &  O_{g}  & 0 \\
0 & 0 & 1
\end{array}
\right) 
\]
for an orthogonal $2\times 2$-matrix $O_{g}$ under a coordinate system. 
 Since $g(c) \subset U_{1}$ from the beginning of (ii),
$g$ does not act on a parallel collection of circles in $A_{1+}$. 
Hence,  $\gamma \ne 0$.

Let $L$ be an annulus bounded by $c$ and $g(c)$ in $A_{1+}$. 
Since there is no bounded component of $L \cap M_{h}$  by Lemma \ref{lem:nobd}, 
$A_{1}$ is the unique component of $M_{h}$.
By Lemma \ref{lem:fullset}, $A_{1+} - A_{1}$ has finitely many components. 
By the action of $\langle g \rangle$, 
the planar surface $A_{1}$ has only one or two ends or infinitely many ends. 
Hence, $A_{1+}= A_{1}$. 
Thus, we conclude $A_{1+} \subset M_{h}$. 

There exists an open neighborhood $N$ of $L$ in $M_h$, and 
$\bigcup_{i\in \bZ} g^i(N) \subset M_h$ covers $A_{1+}$. 
By restricting a Euclidean metric $H^{o}$, we obtain a Euclidean metric on an open set in 
$M_{h}$ containing $\Lambda(R) \cap M_{h}$ and $\bigcup_{i\in \bZ} g^i(L) $.
We obtain a closed set 
\[\Lambda' \subset \bigcup_{i\in \bZ} g^i(N) \cup \Lambda(R) \cap M_{h},\] 
that is foliated by complete $q$-lines and 
\[\Lambda(R) \cap M_{h} \subset \Lambda'\] 
properly.
Also, $\Lambda'$ contains an $\eps$-neighborhood of $\Lambda(R) \cap M_{h}$ in the Euclidean metric. 


The subspace $\Lambda'$ fibers over the surface $\Sigma$ of complete $q$-lines in $\Lambda'$ as before in the beginning of (B).
Then the \hyperlink{term-kcp}{Kuiper completion} $\che \Sigma$  has an affine structure. 
We extend the above fibration 
\[\Pi_{R}: \Lambda_{1}(R) \cap M_{h} \ra S_{R} \hbox{ to }  
\Pi_R: \Lambda' \ra \Sigma\] 
to be denoted by $\Pi_{R}$ again. 
$\Sigma$ properly contains the surface $S_{R}$  bounded by the arc $\alpha$ corresponding to $A_{1+}$.
We take a short geodesic $k$ in $\Sigma$ connecting the endpoints of the short subarc $\alpha_{1}$ in $\alpha$ 
so that they bound a disk in $\Sigma$. $k$ can be extended until it ends in the ideal set of $\che \Sigma$ 
corresponding to complete 
$q$-lines in $\delta_{\infty} \Lambda(R)$. 
We choose a $2$-dimensional crescent $S''$ in $\che \Sigma$ bounded by $k$ containing a $2$-dimensional crescent $S_2$ in $\che S_R$,
and containing $\alpha_{1}$. (See the maximum property in Section 6.2 of \cite{cdcr1}.)
The inverse image $\Pi_R^{-1}(S_2^o)$ has the closure in $\Lambda(R)$ that is a 
bihedral $3$-crescent $S$. 
The union of complete $q$-lines through $S''$ is in $\Lambda'$ if we take $k$ sufficiently close to $\Lambda_{1}(R) \cap M_{h}$. 
The closure of the inverse image $\hat S$ of $S^{\prime \prime o}$ 
in $\che M_{h}$ is a bihedral $3$-crescent in $\Lambda'$ properly containing $S$. 
This contradicts the maximality of $\Lambda(R)$, which is a contradiction to how we defined $\Lambda(R)$ in 
\eqref{eqn:defL}.
\qed

\subsection{The irreducibility of concave affine manifolds.} 


\begin{theorem}\label{thm:concaveirr} 
	Let $N$ be a concave affine $3$-manifold in a connected compact real projective $3$-manifold $M$.
	Suppose that $M$ is neither complete affine nor bihedral. 
	Assume that $M$ has no two-faced submanifold. 
	Then 
	$N$ is irreducible, or $M$ is an affine Hopf $3$-manifold and $M = N$.  
	A \hyperlink{term-camI}{concave affine $3$-manifolds} are always prime. 
\end{theorem} 
\begin{proof} 
	Let $N_h$ denote 
	a component of the inverse image of $N$ in $M_h$. 
	
	Suppose that $N$ is a concave affine $3$-manifold of type I. 
	$N_h$ is a union of an open $3$-hemisphere $R^o$ with a disk $I_{R}^{o}$ in $\partial H$
	where $R$ is a hemispherical $3$-crescent.
	Hence, $N_h$ is irreducible. 
	
	
	Suppose that $N$ is a concave affine $3$-manifold of type II. 
	We assume that $\che M_{h}$ has no hemispherical $3$-crescent. (See Hypothesis \ref{hyp:nohem}.)
	We follow the proof of Theorems \ref{thm:caffI} and \ref{thm:caffII} in 
	Section \ref{subsec:caff}. 
	We divide into cases (A) and (B).
	
	Let $R$ be a bihedral $3$-crescent so that $N_h  \subset \Lambda(R)$.
	\begin{enumerate}
		\item[(A)] Suppose that there are three mutually overlapping bihedral $3$-crescents $R_1, R_2,$ and $R_3$ with 
		$\{I_{R_i}|i=1,2,3\}$ in general position. 
	\end{enumerate}
	We can have two possibilities: 
	\begin{enumerate}
		\item[(i)] Suppose that there exists a pair of \hyperlink{term-opp}{opposite} 
		bihedral $3$-crescents $S_1, S_2 \sim R$. (See Section \ref{subsub:IIAi}.)
		\item[(ii)] There is no such pair of bihedral $3$-crescents. (See Section \ref{subsub:IIAii}.)
	\end{enumerate} 
	
	By Lemma \ref{lem:casei}, (i) implies that $M$ is an affine Hopf manifold. 
	In this case, $\che M$ equals a closed hemisphere and equals $\Lambda(R)$ for 
	a crescent $R$. Thus, $M = N$ for a concave affine manifold $N$. 
	By Proposition \ref{prop:Hopf}, $N$ is a prime $3$-manifold.  
	
	We now work with (ii). 
	By Theorem \ref{thm:bdsph}, $\Bd \Lambda(R) \cap M_{h}$ has no sphere boundary. 
	Any $2$-sphere in $\Lambda(R) \cap M_h$ can be isotopied into $\Lambda_1(R) \cap M_h$. 
	Recall that we showed using Lemma \ref{lem:nobd} 
	in the beginning of (A)(ii) in the proof of Theorem \ref{thm:caffII}
	 that $K$ cannot be bounded in $H^o$ (see Section \ref{subsub:IIAii}).
	We obtain $K \cap \partial H\ne \emp$. 
	Since $H- K$ deformation retracts to $\partial H - K$ by projection 
	from a point of $K^o$, $H - K$ is contractible. 
	Thus, $\Lambda_1(R) \cap M_h$ is contractible and every immersed sphere is null-homotopic.
	
	Now we go to the case (B) in the proof of Theorem \ref{thm:caffII} 
		where $\Lambda(R)$ is a union of the segments whose developing image end at the antipodal pair $q, q_-$ 
		(see Section \ref{subsub:B}).
	Since the interior of $\Lambda(R)\cap N_h$ fibers over an open surface with fiber homeomorphic to real lines, 
	$N$ is irreducible. 
\end{proof}


\subsection{Toral $\pi$-submanifolds} \label{subsec:toralpi}


\begin{lemma} \label{lem:toralpiI} 
A toral $\pi$-submanifold $N$ of type I is homeomorphic to a solid torus or a solid Klein-bottle
and is a concave affine $3$-manifold of type I. 
\end{lemma}
\begin{proof} 
By Theorem \ref{thm:bdsph}, there is no boundary component of $N$ homeomorphic to a sphere or a real projective plane. 

By Definition \ref{defn:toralpi}, $N$  is covered by $R^o \cup I_{R}^o - \{x\}$ for a hemispherical $3$-crescent $R$ and 
$x \in I_R$ and hence is a concave affine $3$-manifold of type I. 

Since the deck transformation group acts on the annulus $I_R^o -\{x\}$ properly discontinuously and freely, 
the deck transformation group of $\tilde N$ is isomorphic to a virtually infinite-cyclic group by Lemma \ref{lem:actann}. 
By Lemma \ref{lem:vic}, we are done. 
\end{proof}

\begin{lemma}\label{lem:toralpi} 
Suppose that $\che M_{h}$ has no hemispherical $3$-crescent. {\rm (}See Hypothesis \ref{hyp:nohem}.{\rm )}
Let $N$ be a concave affine $3$-manifold of type II in $M$ covered by $\Lambda(R) \cap N_{h}$ for 
a bihedral $3$-crescent $R$.  We suppose that
\begin{itemize}
\item The Kuiper completion $\che N_h$ of some cover $N_h$ of the holonomy cover of $N$
contains a \hyperlink{term-tbc}{toral bihedral $3$-crescent} 
$S$ where a deck transformation $g$ acts on $S^o \cup I_S^o -\{x\} \subset N_h$,
fixing a point $x \in I_S^o $ as an attracting fixed point.
\item The deck transformation group of $N$ is virtually infinite-cyclic.
\item $\dev_h|\Lambda(R) \cap N_h: \Lambda(R) \cap N_{h} \ra H^o - K^o$ is an embedding to its image containing $H^o-K$ 
for a compact convex domain $K$ in the $3$-hemisphere $H$
with $K^o \subset H^o$, $K^o \ne \emp$. 
\end{itemize}
Then $N$ contains a unique maximal \hyperlink{term-tsp}{toral $\pi$-submanifold of type II}, 
homeomorphic to a solid torus or a solid Klein-bottle, 
and the interior of every bihedral $3$-crescent in $\che N_h$ meets
the inverse image of the toral $\pi$-submanifold in $N_h$.
\end{lemma}
\begin{proof} 
By Theorem \ref{thm:bdsph}, there is no boundary component of $N$ homeomorphic to a sphere or a real projective plane. 
By definition, $N_h = \Lambda(S) \cap M_{h} $ for a bihedral $3$-crescent $S$. 
We obtain a toral bihedral $3$-crescent $R$ in $\che N_h$. 



By assumption, $\Gamma_N$ is virtually infinite-cyclic.
Two bihedral $3$-crescents $R_1$ and $R_2$ are not \hyperlink{term-opp}{opposite}
since $K^o \ne \emp$ holds. 
Let $R_1$ and $R_2$ be two  \hyperlink{term-tbc}{toral bihedral $3$-crescents} such that $R_{1}, R_{2}\sim R$. Let $R'_i$ denote 
$R_i^o \cup I_{R_i}^o - \{x_i\}$ for a fixed point $x_i$ of the action of an infinite order generating deck transformation 
$g_i$ acting on $R_i$ so that $R'_i/\langle g_i \rangle$ is homeomorphic to a solid torus. 
Here, $g_{i}$ is the deck transformation \hyperlink{term-ass}{associated to} $R_{i}$. 
Let $F_i$, $i=1, 2$, denote the compact fundamental domain of $R_i'$. 
Then the set 
\[G_i:= \{g \in \Gamma_N| g(F_i) \cap F_i \ne \emp\}, i=1,2, \] is finite. 
We can take a finite index normal subgroup $\Gamma'$ of the virtually infinite-cyclic group 
$\Gamma_N$ so that $\Gamma' \cap G_i :=\{e\}$ for both $i$. 
Then a cover of the compact submanifold $R'_i/\langle g_i \rangle \cap \Gamma'$ is embedded in $N_h/\Gamma'$. 
Thus, there is some cover $N_1$ of $N$ so that these lift to embedded submanifolds. 


We denote these in $N_{1}$ by $T_1$ and $T_2$. 
We may assume that $T_i = R'_i/\langle g_i' \rangle$. 
$x_i$ is the {\em fixed point } of $R_i$ and $g_i'$ acts on $R_i'$. 

Suppose that they overlap. Then $R_1 \cap R_2$ is a component of $R_1 - I_{R_2}$  by Theorem 5.4 of \cite{psconv}. 
Considering $T_1 \cap T_2$ that must be a solid torus not homotopic to a point in each $T_i$, 
we obtain that a nonzero power of $g_1'$ and a nonzero power of $g_2'$ are equal.
Therefore, $x_1 = x_2$ and $g_i'$ fixes the point $x_1=x_2$. 

We say that two toral bihedral $3$-crescents $R_1$ and $R_2$ are {\em equivalent} if they overlap in a cover of a solid torus in $N_h$. 
This relation 
generates an equivalence relation of toral bihedral $3$-crescents. We write $R_1 \cong R_2$. 


Let $S$ be a toral bihedral $3$-crescent in $\che N_h$. 
By this condition, 
$x= x_i$ for every fixed point $x_i$ of a toral bihedral $3$-crescent $R_i, R_i \cong S$
where $x$ is fixed by $g_i'$ \hyperlink{term-ass}{associated to} $R_{i}$. 
We define 
\[\hat \Lambda(S):= \bigcup_{R' \cong S} R', \quad \delta_\infty \hat \Lambda(S) := \bigcup_{R'\cong S} \alpha_{R'}.\] 

We claim that $\hat \Lambda(S) \cap N_h$ covers a compact submanifold in $N$: 
Let $T$ be any bihedral $3$-crescent in $\che N_h$ where $g$ acts with $x$ as an attracting fixed point. 
Then $T - \clo(\alpha_T) -\{x\} \subset N_h$ as in cases (II)(A)(ii) or (II)(B)(i).
(See Sections \ref{subsub:IIAii} and \ref{subsub:B}.) 



Since there are no \hyperlink{term-tfsm}{two-faced submanifolds}, we see that 
either 
\[\hat \Lambda(S) = g(\hat \Lambda(S))
\hbox{ or } \hat \Lambda(S) \cap g(\hat \Lambda(S)) \cap N_h = \emp \hbox{ for } g \in \Gamma_N.\]
(This follows as in Lemma 7.2 of \cite{psconv}.)
We can also show that the collection
 \[\{g(\hat \Lambda(S)\cap N_h)| g\in \Gamma_N\}\] is locally finite in $M_h$ as we did for 
$\Lambda(S) \cap M_h$ in Chapter 9 of \cite{psconv}. 
Hence, the image of $\hat \Lambda(S)\cap N_h$ is closed in $M_h$, and it covers a compact submanifold in $M$.



Since $\hat \Lambda(S)$ is a union of segments from $x$ to 
\[\delta_\infty\hat \Lambda(S):=\delta_\infty \Lambda(R) \cap \hat \Lambda(S) ,\] 
$\Bd \hat \Lambda(S) \cap N_h$ is on a union $L$ of such segments from $x$ to $\clo(\delta_\infty \hat \Lambda(S))$ 
passing the set. The open line segments are all in $N_h$ as they are  in toral bihedral $3$-crescents. 
Since $\hat \Lambda(S)$ is canonically defined, the virtually infinite-cyclic group $\Gamma_N$ acts on the set. 
Also, $\hat \Lambda(S) \cap N_h$ is connected since we can apply the above paragraph to $3$-crescents in $\hat \Lambda(S)$ also.

The interior of $\hat \Lambda(S) \cap M_h$ is a union of open segments from $x$ to an open surface $\delta_\infty \hat \Lambda(S)$. 
The surface cannot be a sphere or a real projective plane since a toral $\pi$-submanifold has boundary. 
Since $\delta_\infty \hat \Lambda(S)$ is the complement of $\partial H$ of a compact convex set, 
it is thus homomorphic to a $2$-cell. 
Therefore, the interior of $\hat \Lambda(S)\cap M_h$ is homeomorphic to a $3$-cell. 
We showed that $\hat \Lambda(S)\cap M_h$ covers a compact 
submanifold in $N$. 
We call this $T_N$. 

Clearly, $T_N$ is maximal among the toral $\pi$-submanifolds 
in $N$ since it includes all toral bihedral $3$-crescents. 
If $T_N \subset T'_N$ for any other toral $\pi$-submanifold $T'_N$, 
then the toral bihedral $3$-crescents in a universal cover of $T'_N$ 
overlapping with ones in $\hat \Lambda(S)$ must be in $\hat \Lambda(S)$
since we showed that in the above part of the proof. 
Hence, by considering chains of overlapping bihedral $3$-crescents, we obtain
$T'_N \subset N$ and $T'_N = T_N$. This shows that $T_N$ is a 
maximal toral $\pi$-submanifold. 

Since $N$ has the virtually infinite-cyclic holonomy group, and $\dev_{h}|N_{h}$ is injective, 
we obtain that the holonomy group image of the deck transformation group acting on $\hat \Lambda(S) \cap M_h$ is virtually infinite-cyclic. 
Since the holonomy homomorphism is injective, the deck transformation group acting on $\hat \Lambda(S) \cap M_h$ is virtually infinite-cyclic. 
By Lemma \ref{lem:vic}, the toral $\pi$-submanifold is homeomorphic to a solid torus or a solid Klein bottle.


Now, we go to the final part: 
We assumed that $\dev_h| \Lambda(R)\cap N_h$ is an injective map into the complement of 
a convex domain $K^o$ in $H^o$. 
Thus $\dev_h(\hat \Lambda(S) \cap N_h)$ is the complement of 
a domain $K'$ in $H$ where $K' \supset K^o$ and the closure of $K'$ is convex
and is a union of segments from $\dev_h(x)$ to a domain in $\partial H$. 
Given any bihedral $3$-crescent $R_1$ in $\Lambda(S)$, 
suppose that the open $3$-bihedron $\dev_h(R_1^o)$ does not meet $\dev_h(\hat \Lambda(S))$. 
Then $\dev_h(\alpha_{R_1})$ and $\dev_h(\alpha_{T})$ for a toral bihedral $3$-crescent $T$, 
$T \cong S$, have to be $2$-hemispheres in $\partial H$ antipodal to each other. 
Let $g_T$ denote the deck transformation acting on $T^o \cup I_T^o - \{x\}$ for an attracting fixed point $x$ of $g_T$. 
Then \[g_T^i(R_1) \subset g_T^j(R_1) \hbox{ for } i < j\] by Proposition 3.9 of \cite{psconv} 
since their images overlap and the image of the latter set contains the former one and $\dev_h| N_h$ is injective. 
Hence, the closure of $\bigcup_{i \in \bN} g_T^i(R_1)$ is another toral bihedral $3$-crescent since $g_T$ acts on it. 
Then $T$ and $R$ are \hyperlink{term-opp}{opposite}. 
This is a contradiction since $K$ then has to have an empty interior.
We assumed otherwise in the premise. 

 


\end{proof}

\subsection{Proof of Theorem \ref{thm:main1}}

\begin{proof} 
Given a connected compact real projective $3$-manifold $M$ with empty or convex boundary, 
if $M$ has  a non-$\pi_{1}$-injective component of the two-faced totally geodesic submanifold of type I, then $M$ is an affine Hopf manifold
by Theorem \ref{thm:tfaced}. 

Now suppose that $M$ is not an affine Hopf manifold. 
We \hyperlink{term-spl}{split along} the two-faced totally geodesic submanifolds of type I now to obtain $M^s$.  
Theorem \ref{thm:caffI}  implies the result. 

To complete, we repeat the above argument for  the two-faced totally geodesic submanifold of type II,
and Theorem \ref{thm:caffII} implies the result.
\end{proof}


\section{Toral $\pi$-submanifolds and the decomposition} \label{sec:toraldec}







We now prove a simpler version of Theorem \ref{thm:main3}.
\begin{theorem} \label{thm:decompose}
Suppose that a connected compact real projective $3$-manifold $M$ with empty or convex boundary. 
Suppose that $M$ is neither complete affine nor bihedral, and $M$ is not an affine Hopf $3$-manifold.
Suppose that $M$ has no two-faced submanifold of type I,  and $M$ has no concave affine $3$-manifold of type I
with boundary incompressible to itself. 
Then the following hold\,{\rm :}
\begin{itemize}
\item each concave affine submanifold of type I in $M$ with compressible
boundary contains a unique toral $\pi$-submanifold $T$ of type I where
$T$ has a compressible boundary.
\begin{itemize}
\item There are finitely many disjoint toral $\pi$-submanifolds \[T_1, \dots, T_n\]  obtained by taking one from 
each of the concave affine submanifolds in $M$ with compressible boundary. 
\end{itemize} 
\item We remove $\bigcup_{i=1}^{n} \inte T_{i}$ from $M$. Call $M'$ the resulting real projective manifold with convex boundary. 
Suppose that $M'$ has no two-faced submanifold of type II,  and $M'$ has no concave affine $3$-manifold of type II
with boundary incompressible to itself. 
\begin{itemize} 
 \item  Each concave affine submanifold of type II in $M'$ with compressible
boundary contains a unique toral $\pi$-submanifold $T$ of type II where
$T$ has a compressible boundary.
\item There are finitely many disjoint toral $\pi$-submanifolds \[T_{n+1}, \dots, T_{m+n}\]  obtained by taking one from 
each of the concave affine submanifolds in $M'$ with compressible boundary. 
\end{itemize} 
\item $M - \bigcup_{i=1}^{n+m} \inte T_i$ is $2$-convex.
\end{itemize}
\end{theorem} 
\begin{proof} 
If $N$ is a concave affine $3$-manifold of type I with compressible boundary into $N$, 
then its universal cover is in a hemispherical $3$-crescent, 
and $N$ is homeomorphic to a solid torus and is a toral $\pi$-submanifold by Lemma \ref{lem:toralpiI}. 
These concave affine $3$-manifolds are mutually disjoint. 

We remove these and denote the result by $M'$.
Then $M - \bigcup_{i=1}^n \inte T_i$ has totally geodesic boundary.
The cover $M'_{h}$ of $M'$ is given by removing the inverse images of $T_{1}, \dots, T_{n}$ from $M_{h}$.
We take a Kuiper completion $\che M'_{h}$ of $M'_{h}$. 
Now we consider when $N$ is a concave affine $3$-manifold arising from bihedral $3$-crescents in $\che M'_h$.
We obtain toral $\pi$-submanifold II in $N$ by Lemma \ref{lem:toralpi}.

From $M'$ we remove the union of the interiors of toral $\pi$-submanifolds $T_{n}, \dots, T_{n+m}$.
Then $M - \bigcup_{i=1}^{n+m} \inte T_i$ has a convex boundary as $P_i$ has concave boundary. 

We claim that this manifold $M - \bigcup_{i=1}^{n+m} \inte T_i$ is $2$-convex. Suppose not. Then by Theorem 1.1 of \cite{psconv},  
we obtain again a $3$-crescent $R'$ in the \hyperlink{term-kcp}{Kuiper completion} of 
$M_h - p_{h}^{-1}(\bigcup_{i=1}^{n+m} \inte T_i)$. 
The $3$-crescent $R'$  has 
the interior disjoint from ones we already considered.
However, Theorem \ref{thm:caffII} shows that $R^{\prime o}$ 
must meet the inverse image $p_{h}^{-1}(\bigcup_{i=1}^n \inte T_i)$,
which is a contradiction. 

Lemma \ref{lem:toralpi} shows that each $T_i$ is homeomorphic to a solid torus or a solid Klein bottle. 
\end{proof}

\begin{proof}[Proof of Theorem \ref{thm:main3}.] 
We may assume that $M$ is not complete or bihedral since then $M$ is convex and the conclusions are true. 
As stated, $\che M_{h}$ does not contain any hemispherical $3$-crescent. 
By Theorem \ref{thm:decompose}, $M$ either is an affine Hopf $3$-manifold, or  
$M^{s}$ decomposes into \hyperlink{term-camI}{concave affine $3$-manifolds} with boundary  incompressible  into themselves of type I, 
\hyperlink{term-tsp}{toral $\pi$-submanifolds of type I}, 
and $M^{(1)}$. 

Now $M^{(1)s}$ decomposes into concave affine manifolds of type II with  boundary compressible or incompressible to themselves.  
Theorem 0.1 of \cite{uaf} shows that a $2$-convex affine $3$-manifold is irreducible. 
Toral $\pi$-submanifolds and concave affine $3$-manifolds of type II with incompressible boundary are irreducible or prime
by Lemma \ref{lem:toralpi} and Theorem \ref{thm:concaveirr}. 

\end{proof} 

\begin{proposition} \label{prop:toralemb} 
Let $M$ be a connected compact real projective manifold with convex boundary. 
Suppose that $M$ is not an affine Hopf manifold. 
Then a toral $\pi$-submanifolds of type I in $M^{s}$ is disjoint from the inverse images in $M^{s}$ of 
the two-faced submanifolds in $M$ of type I. Hence, it embeds into $M$. 
Furthermore, 
a toral $\pi$-submanifolds of type II is disjoint from the inverse images in $M^{(1)s}$ of 
the two-faced submanifolds in $M^{(1)}$ of type I. 
And its image in $M^{s}$ is also disjoint from the two-faced submanifolds in $M$ of type I. 
Hence, it embeds into $M$. 
\end{proposition}
\begin{proof} 
Suppose that $M^{s}$ contains a toral $\pi$-submanifold $N$ of type I. 
Then $N^{o}$ embeds into $M$. Then $N^{o}$ is disjoint from the two-faced submanifold $F$ of type I in $M$
by the definition of concave affine manifolds of type I. 
Suppose that the unique boundary component $\partial N$ of $N$ meets the submanifold $F'$ in $M^{s}$ 
mapped to $F$. 
Then since $F'$ is totally geodesic and $\partial N$ is concave, it follows that $\partial N \subset F'$. 
Now, $F$ is non-$\pi_{1}$-injective 
since $\partial N$ is compressible in $N$. 
Theorem \ref{thm:tfaced} shows that $M$ is an affine Hopf $3$-manifold. 
Hence $\partial N$ is disjoint from $F'$, and $N$ embeds into $M$. 

Suppose that $M^{(1)s}$ contains a toral $\pi$-submanifold $N$ of type II. Then $N^{o}$ embeds into 
$M^{(1)}$. Suppose that $\partial N$ meets the submanifold $F'_{2}$ in $M^{(1)s}$ 
mapped to the two-faced submanifold $F_{2}$ of type II in $M^{(1)}$. 
As above, $\partial N \subset F'_{2}$ for the inverse image $F'_{2}$ in
$M^{(1)s}$ of $F_{2}$, and $\partial N$ covers a component $F_{3}$ of $F_{2}$.  
Since $\partial N$ is compressible in $N$, 
Theorem \ref{thm:tfaced} shows that $M$ is an affine Hopf $3$-manifold. 
Hence, $F'_{2}\cap \partial N = \emp$, and 
$N$ embeds into $M^{(1)}$. Call the image by the same name. 

Again $N^{o}$ is disjoint from $F'$. As above $N$ is disjoint from $F'$ or $\partial N \subset F'$. 
In the second case, Theorem \ref{thm:tfaced} shows that $M$ is an affine Hopf $3$-manifold. 
Thus, $N$ embeds into $M$.

\end{proof}

\begin{proof}[Proof of Corollary \ref{cor:main4}.] 

Assume that $M$ is not an affine Hopf $3$-manifold. 
By Proposition \ref{prop:toralemb},  if there exists a toral $\pi$-submanifold in $M^{(1)s}$ or in $M^{s}$, then there is one in $M$. 
Thus, the premise implies that there is no toral $\pi$-submanifold in $M^{s}$ and $M^{(1)s}$. 


Hence, $M^{(1)s}$ decomposes into concave affine $3$-manifolds of type II with incompressible boundary and 
$2$-convex affine $3$-manifolds. Since these are irreducible and each boundary component is not homeomorphic to a sphere
by Theorem \ref{thm:bdsph}, $M^{(1)s}$ is irreducible. 
Since two-faced submanifold $F_{2}$ is $\pi_{1}$-injective by Theorem \ref{thm:tfaced}, 
any sphere $S$ in $M^{(1)}$ meets $F_{2}$ in a disjoint union of circles after perturbations. 
Any disk component of $S - F_{2}$ can be isotopied away since such a disk lifts to one in $M^{(1)s}$
with boundary in  the incompressible surface $F'_{2}$.
By induction, we may assume that $S$ is in 
$M^{(1)} -F_{2}$. Hence, it bounds a $3$-ball. 
Thus, we obtained that $M^{(1)}$ is irreducible as well. 

Now, $M^{s}$ is a union of $M^{(1)}$ and a concave affine manifold of type I with incompressible boundary. 
Similar argument shows that $M^{s}$ and $M$ are irreducible and $M$.

\end{proof}

\begin{proof}[Proof of Corollary \ref{cor:main5}.] 
Suppose that $M$ has an embedded sphere $S$. The domain $\Omega$ contains a lift $S'$ of $S$.
If $S$ is nonseparating, then Corollary \ref{cor:cyclic}  shows that $M$ is an affine Hopf manifold. 

Suppose that $S$ is separating. 
Then $S$ bounds a $3$-ball $B$ in $M$ by Theorem 1.1 of Wu \cite{Wu}. 
\end{proof}

\appendix 

\section{Contraction maps} 

Here, we will discuss contraction maps in $\bR^{n}$. 
A {\em contracting map} $f:X \ra X$ for a metric space $X$ with metric $d$ is a map so that 
$d(f(x), f(y)) < d(x, y)$ for $x, y \in X$. 

\begin{lemma} \label{lem:Katok} 
A linear map $L$ has the property that all the norms of the eigenvalues are $< 1$. 
if and only if $L$ is a contracting map for the distance induced by a norm. 
\end{lemma} 
\begin{proof} 
See Corollary 1.2.3 of Katok \cite{Kt}.
\end{proof}

\begin{proposition} \label{prop:hopf}
$\langle g\rangle $ acts on $\bR^{n} -\{O\}$ {\rm (}reps. $U-\{ O\}$ for the 
upper half space $U \subset \bR^{n}${\rm )}
properly if and only if the all the norms of the eigenvalues of $g$  are $> 1$ or $< 1$. 
\end{proposition} 
\begin{proof}
Suppose that $\langle g\rangle $ acts on $\bR^{n} -\{O\}$ properly. For a sphere $S=\SI^{n-1}$, 
$g^{n}(S)$ is inside a unit ball $B$ for some integer $n$ by the properness of the action. 
This implies that $g^{n}(B) \subset B$, and the norms of the eigenvalues of $g^{n}$ are $<1$. 
The case of the half  space $U$ is similar. 

For the converse, 
by replacing $g$ with $g^{-1}$ if necessary, we assume that all norms of eigenvalues $< 1$. 
Lemma \ref{lem:Katok} implies the result.
\end{proof} 

\begin{proposition}\label{prop:diskf} 
Let $D$ be a domain in $\SI^{n}$. 
Let $g$ be a projective automorphism of $\SI^{n}$  acting on $D$ and
an affine patch $\bR^{n}$.  
We assume the following\,{\rm :} 
\begin{itemize}
\item $S$ is a compact connected subset of $D$ so that $D-S$ has two components $D_{1}$ and $D_{2}$
where $D_{1}$ is bounded in an affine path $\bR^n$ in $\SI^n$. 
\item $g$ acts with a fixed point $x \in \bR^{n}$ in the closure of $D_{1}$. 
\item $g(S) \subset D_{1}$. 
\item Every complete affine line containing $x$ meets $S$ at at least one point.
\item $D_{1} \subset \{x\}\ast S$ where $\{x\}\ast S$ is the union of all segments from $x$ ending at $S$. 
\end{itemize} 
Then $x$ is the global attracting fixed point of $g$ in $\bR^{n}$. 
\end{proposition} 
\begin{proof} 
Choose the coordinate system on $\bR^n$ so that $x$ is the origin. 
Let $L(g)$ denote the linear part of the $g$ in this coordinate system. 
Suppose that there is a norm of the eigenvalue of $L(g)$ greater than or equal to $1$. 
Then there is a real eigensubspace $V$ of dimension $1$ or $2$ associated to 
an eigenvalue of norm $\geq 1$. 
We obtain $S_{V}:= V \cap S \ne \emp$ by the fourth assumption. 
Let $\Theta(S_{V})$ denote the set of directions of $S_{V}$ from $x$. 
$L(g)$ acts on the space of directions from $x$. 
Since $\{x\} \ast g(S) \subset \{x\} \ast S$, we obtain $L(g)(\Theta(S_{V})) \subset \Theta(S_{V})$.
Hence, $\Theta(S_{V})$ is either the set of a point, the set of a pair of antipodal points, or a circle. 
Since $V$ has an invariant metric, there is a point $t$ of $S_{V}$ where a maximal radius of $S_{V}$ takes place. 
Then $g(t) \in g(S_{V})$ must meet $D_{2} \cup S_{V}$, a contradiction.  



Thus, the norms of eigenvalues of $L(g)$ are $< 1$. 
By Lemma \ref{lem:Katok}, $L(g)$ has a fixed point $x$ as an attracting fixed point. 
The conclusion follows. 
\end{proof} 

\section{The boundary of a concave affine manifolds is not strictly concave.} 

The following is the easy generalization of the maximum property in Section 6.2 of \cite{cdcr1}.
A {\em strictly concave point} of a manifold $N$ is a point $y$ where 
 no totally geodesic open disk $D$ containing $y$, $y\in D^{o}$, and 
$D -\{y\} \subset N^{o}$.

\begin{theorem} \label{thm:nscc}
Let $N$ be a concave affine $3$-manifold of type II in $M$. 
Then $\partial N$ has no strictly concave point. 
\end{theorem} 
\begin{proof} 
Let $M_{h}$ be a cover as in the main part of the paper. 
Suppose that there is a disk $D$ as above. 
Then if $y$ is a boundary point of $M_{h}$, then $D$ must be in $\partial M_{h}$ by geometry. 
This contradicts the premise since $D -\{y\} \subset N^{o}$.

Suppose that $N$ is covered by $\Lambda(R) \cap M_{h}$. 
Since $y$ is not a boundary point of $M_{h}$, 
we take a convex compact neighborhood $B(y)$ of the convex point $y$
so that $\dev_{h}(B(y))$ is an $\eps$-$\bdd$-ball for some $\eps > 0$. 
Then $B(y) - \Lambda(R)$ is a properly convex domain with the image
$\dev_{h}(B(y) - \Lambda(R))$ is properly convex. 
For each point $z \in \Bd \Lambda(R)\cap B(y)$, 
let $S_{z}$, $S_{z} \sim R$, be a bihedral $3$-crescent containing $z$. 
Since $\Lambda(R)$ is maximal, 
$\dev_{h}(I_{S_{z}})$ is a supporting plane at $\dev_{h}(z)$ of $\dev_{h}(B(y) - \Lambda(R))$.

We perturb a small convex disk $D \subset I_{S_{y}}$ containing $y$ away from $y$, 
so that the perturbed convex disk $D'$ is such that the closure of $D' \cap B(y) - \Lambda(R)$ is a small compact disk $D''$ 
with 
\[\partial D'' \subset \Bd \Lambda(R) \cap M_{h} \hbox{ and } D^{\prime \prime o} \cap \Lambda(R) = \emp.\] 
Moreover, $\partial D''$ bounds a compact  disk $B'$ in $\Bd \Lambda(R) \cap B(y)$.
Choose a point $z_{0}$ in the interior of $D''$. 
For each point $z \in B'$, $I_{S_z}^o$ is transversal to $\overline{z_0z}$ 
since $z_0 \not\in S_z$. 
Since $S_z^o$ is further away from $z_0$ then $z$, 
we can choose a maximal segment 
$s_{z} \subset S_{z}$ starting from $z_{0}$ passing $z$ 
ending at a point $\delta_{+} s_{z}$ of $\alpha_{S_{z}}$. 
We obtain a compact $3$-ball $B_{z_{0}} = \bigcup_{z\in B'} s_{z}$
with its boundary in $\delta_{\infty }\Lambda(R)$. 
The boundary is the union of $D_{z_{0}} := \bigcup_{z\in \partial D''} s_{z}$, a compact disk, and an open disk 
\[\alpha_{z_{0}} := \bigcup_{z\in B^{o}_{z_{0}}} \delta_{+}s_{z} \subset \delta_{\infty} \Lambda(R).\] 

The image of $\dev_h(\delta_{\infty }\Lambda(R)) \subset \partial H$ 
for a hemisphere $H$ as shown in Page 61 of \cite{psconv}. 
The boundary of $\dev_h(D_{z_0})$ is in $\partial H$. 
For any $\epsilon$, $\epsilon > 0$, 
by taking $D'$ sufficiently close to $D$, we obtain that 
$\dev_{h}(D_{z_{0}})$ is $\epsilon$-$\bdd$-close to $\dev_{h}(I_{S_{y}})$.
It follows that $\dev_{h}(\alpha_{z_{0}})$ 
is $\epsilon$-$\bdd$-close to $\dev_{h}(\alpha_{S_{y}})$.
Hence, $\dev_{h}(B_{z_{0}})$ is a bihedron and $B_{z_{0}}$ is a bihedral $3$-crescent. 

Since $B_{z_{0}}$ is a crescent $\sim S_{y}, S_{y} \sim R$, 
we obtain $B_{z_{0}} \subset \Lambda(R)$. This contradicts our choice of $y$ and $D''$.  
\end{proof}



\bibliographystyle{plain}



\end{document}